\documentclass[10pt]{article}
\usepackage{latexsym,epsf,amssymb,subeqn}
\usepackage{graphicx,mathabx,amsthm,color}
\usepackage{repstyle}

\newcommand{\eq}{\begin{equation}}
\newcommand{\eeq}{\end{equation}}
\newcommand{\R}{\mathbb R}
\newcommand{\C}{\mathbb C}

\newcommand{\Real}{{\rm Re}}

\newcommand{\Oh}{{\cal O}}
\newcommand{\sfrac}[2]{\mbox{\large{$#1\over#2$}}}
\newcommand{\mfrac}[2]{\mbox{\Large{$#1\over#2$}}}
\newcommand{\sDelta}{{\mbox{\footnotesize{$\Delta$}}}}
\newcommand{\dt}{\sDelta t}
\newcommand{\dx}{\sDelta x}
\newcommand{\dy}{\sDelta y}

\newcommand{\eps}{\varepsilon}
\newcommand{\fhi}{\varphi}

\def\imi{\textbf{\hskip1pt i\hskip1pt}}

\begin{document}

\title{On Multistep Stabilizing Correction Splitting Methods with 
Applications to the Heston Model}
\author{
 W.\,Hundsdorfer%
\footnote{CWI, Science Park 123, Amsterdam, The Netherlands.
% \footnote{CWI, Science Park 123, 1098 XG Amsterdam, The Netherlands.
E-mail: willem.hundsdorfer@cwi.nl}
, \quad
K.J.\,in\,'t\,Hout%
\footnote{
Department of Mathematics and Computer Science, University of Antwerp,
Middelheimlaan 1, Antwerp, Belgium.
Email: karel.inthout@uantwerp.be}
}
\date{\today}
\maketitle

\begin{abstract} \noindent
In this note we consider splitting methods based on linear multistep methods 
and stabilizing corrections. 
To enhance the stability of the methods, we employ an idea of Bruno \& Cubillos 
\cite{BrCu16} who combine a high-order extrapolation formula for the explicit 
term with a formula of one order lower for the implicit terms.
Several examples of the obtained multistep stabilizing correction methods are
presented, and results on linear stability and convergence are derived.
The methods are tested in the application to the well-known Heston model 
arising in financial mathematics and are found to be competitive with 
well-established one-step splitting methods from the literature.

\end{abstract}

\section{Introduction} 
\label{Sect:Intro}

In this note we will discuss a class of splitting methods for solving
initial value problems for ordinary differential equations (ODEs)
\eq
\label{eq:ODE}
u'(t) = F(t,u(t)) \,, \qquad u(0) = u_0 \,,
\eeq
with given $u_0\in\R^M$, $F:\R\times\R^M\rightarrow \R^M$ and
dimension $M\ge1$.
For many problems occurring in practice there is a natural decomposition
\eq
\label{eq:Deco}
F(t,u) \,=\, F_0(t,u) + F_1(t,u) + \cdots + F_s(t,u)
\eeq
in which the separate component functions $F_j$ are more simple than
the whole $F$, and where $F_0$ is a non-stiff or mildly stiff term
that can be treated explicitly in a time stepping method. 
For such problems we will study a class of multistep splitting 
methods with stabilizing corrections, where explicit predictions 
are followed by corrections that are implicit in one of the $F_j$ terms, 
$j=1,2,\ldots,s$.

\subsection{Linear multistep methods with stabilizing corrections}

The splitting methods to be considered produce approximations 
$u_n \approx u(t_n)$ at the step points $t_n = n \dt$.  The 
methods are based on pairs of linear multistep methods: an implicit method
\eq
\label{eq:LMimpl}
u_n \,=\, \sum_{i=1}^k a_i u_{n-i} + \dt \sum_{i=0}^k b_i F(t_{n-i},u_{n-i}) 
\eeq
and an explicit one
\eq
\label{eq:LMexpl}
u_n \,=\, \sum_{i=1}^k \Big( a_i u_{n-i}
+ \dt\, \widehat{b}_i F(t_{n-i},u_{n-i}) \Big) \,,
\eeq
having the same coefficients $a_1, a_2,\ldots, a_k$ and the same order $p$. 

Starting with an implicit method of order $k$, the matching explicit formula 
of the same order can be obtained by extrapolation, replacing the implicit term 
$F(t_n,u_n)$ by a linear combination of $F(t_{n-i},u_{n-i})$ terms, 
$1\le i\le k$. More specifically, with 
$\widehat{c}_i = \prod_{j\neq i}\frac{j}{j-i}$, $1\le i,j\le k$, 
it can be seen from the Lagrange interpolation formula that 
\eq
\label{eq:Lagr1}
\fhi(t_n) = \sum_{i=1}^k \widehat{c}_i \fhi(t_{n-i}) + 
\widehat{C} \dt^k \fhi^{(k)}(t_n)
+ \Oh(\dt^{k+1})
\eeq
for any smooth function $\fhi$, with error constant $\widehat{C}\in\R$. 
Using this extrapolation procedure to replace the implicit term in 
(\ref{eq:LMimpl}) yields the following coefficients for the explicit method:
\eq
\label{eq:Lagr1'}
\widehat{b}_i\,=\, b_i + b_0 \widehat{c}_i \qquad (i=1,2,\ldots,k) \,.
\eeq 

Such a pair of explicit and implicit methods can now be combined to form
a splitting method for problems with decomposition (\ref{eq:Deco}) by 
using the idea of stabilizing corrections. In each step, first a prediction
is made with the explicit method, followed by corrections for the implicit
function components $F_j$ with $j =1,2,\ldots,s$:
\eq
\label{eq:SCLM}
\setlength{\arraycolsep}{1mm}
\left\{
\begin{array}{ccl}
v_{0} &=&
\displaystyle
\sum_{i=1}^k \Big( a_i u_{n-i} + \dt\, \widehat{b}_i F(t_{n-i},u_{n-i}) \Big) \,,
\\[-2mm]
v_{j} &=&
\displaystyle
v_{j-1} + \dt \sum_{i=1}^k (b_i - \widehat b_i)  F_j(t_{n-i}, u_{n-i})
+ \dt\, b_0 F_j(t_n,v_{j}) \,,
\\[2mm]
u_{n} &=& v_{s} \,.
\end{array}
\right.
\eeq
All internal vectors $v_0,v_1,\ldots.v_s$ that appear in this step are 
consistent approximations to the exact solution value $u(t_n)$.
This property ensures that steady-state solutions are maintained by the
scheme, that is, if $F(u_*)=0$ and $u_{n-i}=u_*$ for $i=1,\ldots,k$, then
$u_n=u_*$.
Using the terminology of \cite{Ma90}, this splitting method 
will be called a multistep stabilizing correction method.

For linear problems without explicit terms, the above formula (\ref{eq:SCLM})
is closely related to a class of methods introduced by Douglas \& Gunn 
\cite{DoGu64}. In that paper it was noted that stability properties may 
be improved by using extrapolation of lower order. However, with an
explicit term this will generally lead to a lower order of convergence
of the splitting method. 

To overcome this, we will follow an idea of Bruno \& Cubillos \cite{BrCu16} 
who studied BDF splitting methods for linearized Navier-Stokes equations 
using two different extrapolation formulas in the prediction stage: high-order 
extrapolation for the explicit term $F_0$ and a formula of one order lower 
for the implicit terms $F_1,\ldots,F_s$. This lower order formula will be 
\eq
\label{eq:Lagr2}
\fhi(t_n) = \sum_{i=1}^k \widecheck{c}_i \fhi(t_{n-i}) 
+ \widecheck{C} \dt^{k-1} \fhi^{(k-1)}(t_n)
+ \Oh(\dt^{k}) \,,
\eeq
with coefficients $\widecheck{c}_i$ and error constant $\widecheck{C}$.
In our examples we will take $\widecheck{c}_k=0$ and use Lagrange 
interpolation through the data values $\fhi(t_{n-1}),\ldots,\fhi(t_{n-k+1})$.
Analogous to (\ref{eq:Lagr1'}), let 
\eq
\label{eq:Lagr2'}
\widecheck{b}_i\,=\, b_i + b_0 \widecheck{c}_i \qquad (i=1,2,\ldots,k) \,.
\eeq

With this lower order extrapolation in the prediction stage for the
implicit terms $F_1,\ldots,F_s$ we get a multistep splitting method
of the following form:
\eq
\label{eq:SCLMmod}
\setlength{\arraycolsep}{1mm}
\left\{
\begin{array}{ccl}
v_{0} &=&
\displaystyle
\sum_{i=1}^k \Big( a_i u_{n-i} + \dt\, \widehat{b}_i F_0(t_{n-i},u_{n-i}) 
+ \dt\, \widecheck{b}_i \sum_{j=1}^s F_j(t_{n-i},u_{n-i}) \Big) \,,
\\[-2mm]
v_{j} &=&
\displaystyle
v_{j-1} + \dt \sum_{i=1}^k (b_i - \widecheck{b}_i)  F_j(t_{n-i}, u_{n-i})
+ \dt\, b_0 F_j(t_n,v_{j}) \,,
\\[2mm]
u_{n} &=& v_{s} \,,
\end{array}
\right.
\eeq
again with index $j =1,2,\ldots,s$ in the correction steps.
These correction steps now not only serve to provide stability but also 
the accuracy of the prediction step needs to be improved. We will
refer to these methods as modified stabilizing correction methods.  

For the special case $s=1$, where we have only one implicit term, 
both formulas (\ref{eq:SCLM}) and (\ref{eq:SCLMmod}) reduce to
\eq
\label{eq:IMEX}
u_n \,=\,
\sum_{i=1}^k \Big( a_i u_{n-i} + \dt\, \widehat{b}_i F_0(t_{n-i},u_{n-i}) \Big)
+ \dt\, \sum_{i=0}^k b_i F_1(t_{n-i},u_{n-i}) \,. 
\eeq
This gives the well-known class of implicit-explicit (IMEX) linear multistep 
methods. These methods were originally introduced in \cite{Cr80,Va80}, and
a number of interesting examples can be found in \cite{ARW95}.
In many applications, splittings with more implicit terms appear, in which
case (\ref{eq:SCLM}) and (\ref{eq:SCLMmod}) provide natural generalizations
of these IMEX methods.

\subsection{Outline}

In this paper we will discuss the stabilizing correction multistep 
methods with application to problems arising in financial option valuation. 
First we will present in Section~\ref{Sect:Exas} examples of suitable 
pairs of linear multistep methods. 
The accuracy of the stabilizing correction multistep methods is analyzed 
in Section~\ref{Sect:Accu} for linear problems. It will be seen that 
the methods (\ref{eq:SCLMmod}) may show a local order reduction due 
to stiffness, but it will also be seen that under mild assumptions 
the global errors will not be affected by such order reduction.
Section~\ref{Sect:Stab} contains stability results for 2D parabolic 
problems with cross-derivatives, using an ADI type splitting together 
with explicit treatment of the cross-derivative term. 
In Section~\ref{Sect:Heston} numerical results are presented and 
discussed for a well-known 2D problem from financial mathematics, 
the so-called Heston model for option valuation.
Section~\ref{Sect:Concl} contains some final remarks and conclusions.

\section{Examples}
\label{Sect:Exas}

The following examples fit in the framework outlined in the previous 
section, with coefficients in the prediction stage obtained by extrapolation.
The examples are described by specifying $k$ and the coefficients
$a = (a_1,\ldots,a_k)$, $b = (b_1,\ldots,b_k)$, $\theta = b_0$
of the implicit method together with the coefficients
$\widehat{b} = (\widehat{b}_1,\ldots,\widehat{b}_k)$ and
$\widecheck{b} = (\widecheck{b}_1,\ldots,\widecheck{b}_k)$
for the explicit prediction stage.

\bigskip\noindent
{\it The Douglas method\/}:
The most simple example is obtained for the one-step case,  $k=1$.
With the implicit $\theta$-method, extrapolation yields the forward 
Euler method. The combination with stabilizing corrections is the 
Douglas method, which can be written as (\ref{eq:SCLM}) with 
\eq
\label{eq:Doug}
a = 1 \,, \qquad
b = 1-\theta \,, \qquad
\widehat{b} = 1 \,,
\eeq
and the free parameter $\theta\in[\frac{1}{2},1]$. Originally \cite{Do62}
the method was intended for linear parabolic equations without explicit
term $F_0$. This method has been studied in a number of publications, e.g.\ 
\cite{HoWe07,Hu02,HuVe03}.

\bigskip\noindent
{\it The combination CNLF\/}:
A popular combination of implicit and explicit two-step methods
is found with the implicit trapezoidal rule, written in two-step form, 
and the explicit midpoint method. This leads to (\ref{eq:SCLM}) with
$k=2$ and
\begin{subequations}
\label{eq:CNLF}
\eq
\label{eq:CNLFa}
a = (0, \, 1) \,, \qquad
b = (0, \, 1) \,, \qquad
\widehat{b} = (2, \, 0) \,, \qquad
\theta = 1 \,.
\eeq
The implicit trapezoidal rule and explicit midpoint method are often called
the Crank-Nicolson (CN) method and Leap-frog (LF) method in PDE applications.
For the modification (\ref{eq:SCLMmod}) with lower order extrapolation we get
\eq
\label{eq:CNLFb}
\widecheck{b} = (1, \, 1) \,.
\eeq
\end{subequations}

\smallskip \noindent
{\it BDF2 combinations\/}: 
We will consider the class of implicit second-order two-step methods with
free parameter $\theta = b_0 > 0$ and
\begin{subequations}
\label{eq:BDF2}
\eq
\label{eq:BDF2a}
a = \big(\sfrac{4}{3}, \, -\sfrac{1}{3}\big) \,, \qquad
b = \big(\sfrac{4}{3} - 2\theta, \, -\sfrac{2}{3}+\theta\big) \,.
\eeq
For $\theta=\frac{2}{3}$ this is the familiar BDF2 method. With free
parameter $\theta$ it will be referred to as generalized BDF2, even
though the backward differentiation idea is not so prominent anymore
if $\theta \neq \frac{2}{3}$.
These implicit methods are $A$-stable for $\theta\ge\frac{1}{2}$.
With linear and constant extrapolation we get the coefficients
\eq
\label{eq:BDF2b}
\widehat{b} = \big(\sfrac{4}{3}, \, -\sfrac{2}{3}\big) \,, \qquad
\widecheck{b} = \big(\sfrac{4}{3} - \theta, \, -\sfrac{2}{3}+\theta\big) \,.
\eeq
\end{subequations}

\smallskip \noindent
{\it Adams2 combinations\/}: 
As a further example we consider the class of implicit second-order 
two-step Adams-type methods with free parameter $\theta = b_0 > 0$ and
\begin{subequations}
\label{eq:Adams2}
\eq
\label{eq:Adams2a}
a = \big(1, \, 0\big) \,, \qquad
b = \big(\sfrac{3}{2} - 2\theta, \, -\sfrac{1}{2}+\theta\big) \,.
\eeq
These methods are $A$-stable for $\theta\ge\frac{1}{2}$. Linear and
constant extrapolation gives
\eq
\label{eq:Adams2b}
\widehat{b} = \big(\sfrac{3}{2}, \, -\sfrac{1}{2}\big) \,, \qquad
\widecheck{b} = \big(\sfrac{3}{2} - \theta, \, -\sfrac{1}{2}+\theta\big) \,.
\eeq
The IMEX method (\ref{eq:IMEX}) with $\theta = \frac{1}{2}$ is often 
referred to as CNAB because it combines the explicit Adams-Bashforth 
method with the implicit trapezoidal rule (Crank-Nicolson).
Larger values of $\theta$ have been considered in \cite{ARW95,NeLi79},
see also \cite[p.\,388]{HuVe03}.
\end{subequations}

\bigskip \noindent
{\it BDF3 combination\/}:
In the following we will mainly consider two-step methods with order 
two. Higher orders can be obtained with $k>2$.
As an example of a method with $k=3$ we consider the implicit BDF3 method
with coefficients
\begin{subequations}
\label{eq:BDF3}
\eq
a = \big(\sfrac{18}{11},\, -\sfrac{9}{11},\, \sfrac{2}{11}\big) \,, \qquad
b = \big(0,\, 0,\, 0\big) \,\qquad
\theta = \sfrac{6}{11} \,.
\eeq
The use of quadratic and linear extrapolation leads to
\eq
\widehat{b} = \big(\sfrac{18}{11},\, -\sfrac{18}{11},\, \sfrac{6}{11}\big) \,, 
\qquad 
\widecheck{b} = \big(\sfrac{12}{11},\, -\sfrac{6}{11},\, 0\big) \,.
\eeq
\end{subequations}

\section{Discretization errors and convergence}
\label{Sect:Accu}

The accuracy analysis of splitting methods for stiff ODEs and
semi-discrete systems obtained from partial differential equations
(PDEs) should take stiffness into account.

Let
\eq
\label{eq:fhi_j}
\fhi_j(t) = F_j(t,u(t))
\qquad (j=1,2,\dots,s) \,.
\eeq
In the following it will be assumed that the exact solution $u$ and
the functions $\fhi_j$ are all sufficiently smooth on the time interval
$[0,T]$, with derivatives that are bounded uniformly in the stiffness.  
In general, for problems (\ref{eq:ODE}) with a smooth solution, this 
condition on the functions $\fhi_j$  will hold for suitable splittings.

Further it will be assumed that the implicit linear multistep method is of
order $k$, and the extrapolation procedures satisfy (\ref{eq:Lagr1}),
(\ref{eq:Lagr2}).  The question is whether the multistep splitting
method (\ref{eq:SCLMmod}) will be convergent of order $k$ for stiff
problems, and in particular for problems obtained by spatial discretization
of a PDE. For such stiff problems we will denote by $\Oh(\dt^m)$ a
vector whose norm is bounded by $C\dt^m$ with a constant $C$ independent
of the mesh-width $h$ in the spatial discretization. In the same
fashion, $\Oh(1)$ indicates the norm is bounded uniformly in $h$.

For the error analysis we will consider a step (\ref{eq:SCLMmod}), but now
starting from perturbed values $\widetilde{u}_{n-i}$ together with
perturbations $\rho_j$ in the stages, leading to
\eq
\label{eq:SCLMmod'}
\setlength{\arraycolsep}{1mm}
\left\{
\begin{array}{ccl}
\widetilde{v}_{0} &=&
\displaystyle
\sum_{i=1}^k \Big( a_i \widetilde{u}_{n-i}
+ \dt\, \widehat{b}_i F_0(t_{n-i},\widetilde{u}_{n-i})
+ \dt\, \widecheck{b}_i \sum_{j=1}^s F_j(t_{n-i},\widetilde{u}_{n-i}) \Big)
+ \rho_0 \,,
\\[-2mm]
\widetilde{v}_{j} &=&
\displaystyle
\widetilde{v}_{j-1} + \dt \sum_{i=1}^k
(b_i - \widecheck{b}_i)  F_j(t_{n-i}, \widetilde{u}_{n-i})
+ \dt\, b_0 F_j(t_n,\widetilde{v}_{j}) + \rho_j \,,
\\[2mm]
\widetilde{u}_{n} &=& \widetilde{v}_{s} \,.
\end{array}
\right.
\eeq
As for the $\widetilde{v}_{j}$, also the $\rho_j$ depend on $n$.
If we insert exact solution values for $\widetilde{u}_{n-i}$ and
$\widetilde{v}_j$ then the $\rho_j$ become truncation errors for the
stages.

\subsection{Error recursions}

For the analysis%
\footnote{
The derivation of the error recursions and the formulas for the
stage truncation errors can be done for nonlinear problems, but to
get bounds for the local and global errors many technical assumptions
would be needed in the nonlinear case.
}
it will be assumed that the problem is linear,
\eq
\label{eq:ODElin}
F_j(t,u) = A_j u + g_j(t) \qquad (j=0,1,\ldots,s) \,.
\eeq
The source terms $g_j$ may contain inhomogeneous boundary values of the
underlying PDE. The matrices $A_j$ may contain negative powers the
mesh-width in space $h$, and the same applies to $g_j(t)$ if
inhomogeneous boundary values are included in that term.

Further we use the notations
\eq
\eps_{n-i} = \widetilde{u}_{n-i} - u_{n-i} \,, \qquad
\nu_j = \widetilde{v}_j - v_j \,,
\eeq
and
\eq
Z_j = \dt A_j \,, \qquad Q_j = I - b_0 Z_j \, \qquad
P = Q_1 Q_{2} \cdots Q_s \,.
\eeq

With these notations, subtracting the unperturbed scheme from the
perturbed one leads to the relations
\eq
\setlength{\arraycolsep}{1mm}
\left\{
\begin{array}{ccl}
\nu_0 &=& \displaystyle
\sum_{i=1}^k \Big( a_i \eps_{n-i} + \widehat{b}_i Z_0 \eps_{n-i}
+ \widecheck{b}_i \sum_{j=1}^s Z_j \eps_{n-i}\Big) + \rho_0
\\[-2mm]
\nu_j &=& \displaystyle
\nu_{j-1} + \sum_{i=1}^k (b_i-\widecheck{b}_i) Z_j
\eps_{n-i} + b_0 Z_j \nu_j + \rho_j
\\[4mm]
\eps_n &=& \nu_s \,.
\end{array}
\right.
\eeq
Setting $\sigma_n = \sum_{i=1}^k (b_i-\widecheck{b}_i) \eps_{n-i}$ gives
$\nu_j = Q_j^{-1}\big(\nu_{j-1} + Z_j \sigma_n + \rho_j)$.
It follows that
\eq
\eps_n \,=\, P^{-1} \nu_0
+ \sum_{j=1}^s Q_s^{-1}\!\ldots Q_j^{-1} \big(Z_j \sigma_n + \rho_j\big) \,.
\eeq
Substitution of the expressions for $\nu_0$ and $\sigma_n$ now leads to the
recursion
\eq
\label{eq:GlobErr}
\eps_n \,=\, \sum_{i=1}^k R_i\,\eps_{n-i} + \delta_n
\eeq
with error per step
\eq
\label{eq:LocErr}
\delta_n \,=\, P^{-1} \rho_0
+ \sum_{j=1}^s Q_s^{-1}\!\ldots Q_j^{-1}  \rho_j
\eeq
and matrices
\eq
\label{eq:R_i'}
R_i \,=\, P^{-1}\Big( a_i I + \widehat{b}_i Z_0
+ \widecheck{b}_i \sum_{j=1}^s Z_j \Big)
+ \big(b_i - \widecheck{b}_i\big) \sum_{j=1}^s Q_s^{-1}\!\ldots Q_j^{-1} Z_j \,.
\eeq

These matrices can be written is a more simple form.
By induction with respect to $s$ it can be shown that
\eq
b_0 \sum_{j=1}^s Q_1 Q_2 \cdots Q_{j-1} Z_j \,=\, I - P \,.
\eeq
From this relation it follows that
\eq
\label{eq:R_i}
R_i \,=\, P^{-1}\Big( a_i I + \widehat{b}_i Z_0
+ \widecheck{b}_i \sum_{j=1}^s Z_j +
\mfrac{1}{b_0}\big(b_i - \widecheck{b}_i\big) \big(I - P\big) \Big) \,.
\eeq

In this section we will use the above formulas with
$\widetilde{u}_{n} = u(t_{n})$ for all $n$, so that (\ref{eq:GlobErr}) 
becomes a recursion for the global discretization errors 
$\eps_{n} = u(t_{n}) - u_{n}$.
Then $\delta_n$ will be a local discretization error, introduced in the step
from $t_{n-1}$ to $t_n$. The choice of the vectors $\widetilde{v}_j$ is
free, but it is convenient to take $\widetilde{v}_j = u(t_n)$ to obtain
simple expressions for the residuals $\rho_j$.

\subsection{Stability}

In the following it will be assumed that the space $\R^M$ is equipped with
a suitable norm, and that we have in the induced matrix norm
\eq
\label{eq:BoundQinv}
\|Q_j^{-1}\| \,\le\, \kappa
\quad \mbox{for $j=1,\ldots,s$},
\eeq
with a moderately sized constant $\kappa\ge1$. In many instances this will
hold with $\kappa=1$.

Further it will be assumed that the recursion
\eq
\label{eq:GlobErr0}
\eps_n \,=\, \sum_{i=1}^k R_i\,\eps_{n-i}  \qquad (n\ge k)
\eeq
is stable, in the sense that there is a constant $K\ge1$, not affected by
stiffness, such that
\eq
\label{eq:stability0}
\|\eps_n\| \,\le\, K \max_{0\le i\le k-1}\|\eps_i\|
\eeq
for all $n\ge k$ and arbitrary starting errors $\eps_0,\ldots,\eps_{k-1}
\in\R^M$.

For the recursion (\ref{eq:GlobErr}) with local errors $\delta_n$
this will imply
\eq
\label{eq:stability1}
\|\eps_n\| \,\le\, K\Big( \max_{0\le i\le k-1}\|\eps_i\| +
\sum_{j=k}^n \|\delta_j\| \Big) \,,
\eeq
as can be seen by writing the multistep recursion in a one-step form
in a higher dimensional space, see \cite[p.\,183]{HuVe03}, for example.
If $\|\eps_i\|\le C_0\dt^p$ ($0\le i\le k$) and $\|\delta_j\|\le C_1\dt^{p+1}$
($k\le j\le n$) we now get
$ \|\eps_n\| \le (K C_0 + t_n C_1)\dt^p$,
which is the standard way to demonstrate convergence on finite time
intervals, $t_n\in[0,T]$. As we will see, this convergence argument
will need some refinement for the splitting methods (\ref{eq:SCLMmod})
applied to stiff problems.

Verification of the stability condition (\ref{eq:stability0}) can be
quite difficult in practical situations. In Section~\ref{Sect:Stab} this
condition will be studied for a class of parabolic problems with mixed
derivatives in a von Neumann analysis. Then stability in the discrete
$L_2$-norm follows from the scalar case with $z_0,z_1,\ldots,z_k\in\C$
replacing the matrices $Z_0,Z_0,\ldots,Z_k$. For this scalar case we get the
recursion
\eq
\label{eq:RecScal}
\eps_n \,=\, \sum_{i=1}^k r_i \,\eps_{n-i}  \qquad (n\ge k)
\eeq
with
\eq
r_i \,=\, \mfrac{1}{p} \Big(a_i + \widehat{b}_i z_0 + \widecheck{b}_i
\sum_{j=1}^s z_j + \mfrac{1}{b_0}(b_i-\widecheck{b}_i) (1 - p) \Big) \,,
\quad
p \,=\, \prod_{j=1}^s (1 - b_0 z_j) \,.
\eeq
Stability for this recursion is determined by the roots of the
characteristic polynomial
\eq
\label{eq:CharPol}
\pi(\zeta) \,=\, \zeta^k - \sum_{i=1}^k r_i \zeta^{k-i} \,.
\eeq
The coefficients $r_i$ in this polynomial depend on $z_0,z_1,\ldots,z_s$,
so the same holds for its roots $\zeta_l = \zeta_l(z_0,z_1,\ldots,z_s)$.
The recursion is stable iff this polynomial $\pi$ satisfies the well-known
root condition: all roots $\zeta_l$ have modulus at most one, and those
with modulus one are simple.

\subsection{Stage truncation errors}

The internal vectors $v_0,v_1,\ldots.v_s$ in the step (\ref{eq:SCLMmod})
are all consistent approximations to $u(t_n)$.  Insertion of the exact
solution values $\widetilde{v}_{j} = u(t_n)$ and
$\widetilde{u}_{n-i} = u(t_{n-i})$
into (\ref{eq:SCLMmod'}) yields truncation errors
$\rho_0,\rho_1,\ldots,\rho_s$ in the stages given by
\begin{subequations}
\label{eq:SC_res}
\eq
\label{eq:SC_resA}
\rho_0 \,=\, u(t_n) - \sum_{i=1}^k \Big( a_i u(t_{n-i}) 
+ \dt\, \widehat{b}_i \fhi_0(t_{n-i})
+ \dt\,\widecheck{b}_i \sum_{j=1}^s \fhi_j(t_{n-i})\Big)
\,,
\eeq
\eq
\label{eq:SC_resB}
\rho_j \,=\,
-\dt \sum_{i=1}^k (b_i - \widecheck b_i)  \fhi_j(t_{n-i})
- \dt\, b_0 \fhi_j(t_n)  \qquad (j=1,2,\ldots,s) \,.
\eeq
\end{subequations}

To estimate $\rho_0$ we write
\eq
\rho_0 \,=\, \rho_0^{\rm impl} + \rho_0^{\rm extr}
\eeq
with
\begin{subequations}
\eq
\rho_0^{\rm impl} \,=\,
u(t_n) - \sum_{i=1}^k a_i u(t_{n-i})
- \dt \sum_{i=0}^k b_i u'(t_{n-i}) \,,
\eeq
\eq
\rho_0^{\rm extr} \,=\,
\dt \, b_0 u'(t_n) - \dt \, b_0 \sum_{i=1}^k
\Big( \widehat{c}_i \fhi_0(t_{n-i}) 
+ \widecheck{c}_i \sum_{j=1}^s \fhi_j(t_{n-i}) \Big) \,.
\eeq
\end{subequations}
Since $\rho_0^{\rm impl}$ is the truncation error of the implicit method,
we have $\rho_0^{\rm impl} = \Oh(\dt^{k+1})$.
The term $\rho_0^{\rm extr}$ is due to extrapolation. From
(\ref{eq:Lagr1}), (\ref{eq:Lagr2}), with the error constants $\widehat{C},
\widecheck{C}$, it is seen that
\eq
\rho_0^{\rm extr} \,=\, \dt^{k+1} b_0 \widehat{C} \fhi_0^{(k)}(t_n) 
+ \Oh(\dt^{k+2}) + \dt^k  b_0 \widecheck{C} \sum_{j=1}^s \fhi_j^{(k-1)}(t_n)
+ \Oh(\dt^{k+1}) \,.
\eeq
Consequently we have
\begin{subequations}
\label{eq:ResBounds}
\eq
\rho_0 \,=\, \dt^k  \beta_0 \sum_{j=1}^s \fhi_j^{(k-1)}(t_n)
+ \Oh(\dt^{k+1})
\eeq
with $\beta_0 = b_0 \widecheck{C}$.
From (\ref{eq:Lagr2}) it also follows that
\eq
\rho_j \,=\, -\dt^k  \beta_0 \fhi_j^{(k-1)}(t_n) + \Oh(\dt^{k+1}) \,.
\eeq
\end{subequations}

\subsection{Local error bounds}

Combining (\ref{eq:ResBounds}) with (\ref{eq:LocErr}) gives the
following expression for the local discretization error:
\eq
\label{eq:LocErr'}
\delta_n \,=\, \dt^k \beta_0 \sum_{j=1}^s P^{-1}
\big(I  - Q_1\ldots Q_{j-1} \big) \fhi_j^{(k-1)}(t_n) + \Oh(\dt^{k+1}) \,.
\eeq

For the {\it non-stiff case} we have $Z_j = \Oh(\dt)$ and $Q_j = I + \Oh(\dt)$,
and then it follows that $\delta_n = \Oh(\dt^{k+1})$, which is the
usual local error estimate for a method of order $k$.
However, if the ODE system is stiff, for example if the system is
obtained by spatial discretization of a PDE, then these local errors
need more careful examination.

First of all, let us remark that for the case $s=1$ we get $\delta_n =
\Oh(\dt^{k+1})$, and in this remainder term only derivatives of the
$\fhi_j(t) = F_j(t,u(t))$ are involved.  Therefore, if we have fixed
bounds, not affected by stiffness, for the norms of these derivatives,
then also the local truncation errors will not be affected by stiffness
if $s=1$.

A local error bound $\delta_n = \Oh(\dt^{k+1})$ is also valid for the
methods (\ref{eq:SCLM})  where only high-order extrapolation is used,
that is, $\widecheck{b}_i = \widehat{b}_i$ such that (\ref{eq:Lagr1}) holds.
For this case we can use the above formulas with $D=0$, giving $\beta_0=0$.

However, for the methods (\ref{eq:SCLMmod}) with lower-order extrapolation
such $\Oh(\dt^{k+1})$ local error bounds need no longer to be valid in
general.  For example, if $s=2$ then we have
\eq
\label{eq:LocErr''}
\begin{array}{c}
\delta_n \,=\,
\dt^k \beta_0 \,Q_2^{-1} Q_1^{-1} (I - Q_1)  \fhi_2^{(k-1)}(t_n)
+ \Oh(\dt^{k+1})
\\[2mm]
=\, \dt^k \beta_0 \, (I - b_0Z_2)^{-1} (I - b_0Z_1)^{-1}
b_0 Z_1 \, \fhi_2^{(k-1)}(t_n) + \Oh(\dt^{k+1}) \,.
\end{array}
\eeq
Using (\ref{eq:BoundQinv}) it follows that $\|\delta_n\| = \Oh(\dt^k)$,
but the classical bound $\|\delta_n\| = \Oh(\dt^{k+1})$ will not hold in
general for stiff systems.

\begin{Exa} \rm
Consider the model problem consisting of the 2D heat equation
$$
u_t \,=\, u_{xx} + u_{yy} + f(x,y,t)
$$
on the unit square $\Omega = [0,1]^2$ and $t\ge0$, with given initial condition at time $t=0$ and Dirichlet boundary conditions
$$
u(x,y,t) = \gamma(x,y,t)  \qquad
\mbox{for}\; t\ge0\,, \; (x,y)\in \Gamma = \partial \Omega \,.
$$
Standard discretization on a uniform Cartesian grid with mesh-width $h$
in both directions, $\dx=\dy=h$, leads to a semi-discrete ODE system
$$
u'(t) =  \sum_{j=1}^2 \big(A_j u(t) + g_j(t)\big) + g_0(t) \,,
$$
where $g_0(t)$ is the restriction of the source term $f(x,y,t)$ to the
spatial grid, $A_1\approx \frac{\partial^2}{\partial x^2}$ and $g_1(t)$
contains the boundary data for $x=0$ and $x=1$, and likewise in the
$y$-direction for $A_2\approx \frac{\partial^2}{\partial y^2}$ and $g_2(t)$.

If $\fhi_2^{(k-1)}(t_n) = \upsilon_h$ where $\upsilon_h$ is the restriction to
the grid of a smooth function $\upsilon(x,y)$ that is not equal to zero
at the boundaries $x=0$ or $x=1$, then it can be observed in experiments that
$$
\|\delta_n\|_2\sim \dt^{k+1/4} \,, \qquad
\|\delta_n\|_\infty \sim \dt^k \,, \qquad
$$
in the discrete $L_2$-norm $\|v\|_2 = (\frac{1}{M}\sum_{i=1}^M |v_i|^2)^{1/2}$
and the maximum-norm $\|v\|_\infty = \max_{1\le i\le M} |v_i|$,
respectively, for $v=(v_i)\in\R^M$.
In fact, from a spectral analysis, as in \cite[pp.\,296--300]{HuVe03}, it can
be shown that $\|\delta_n\|_2\sim \log(\dt)\cdot\dt^{k+1/4}$, but
such a logarithmic term is hardly observable in experiments.
This order reduction is caused by the boundary conditions, not by lack 
of smoothness of the solution.
\hfill $\Diamond$
\end{Exa}

\subsection{Global error bounds}

Even in situations where $\delta_n \neq \Oh(\dt^{k+1})$ we can have
$\eps_n =  \Oh(\dt^k)$ uniformly in $n$, that is, convergence of order $k$,
due to damping and cancellation effects. This happens for many one-step
splitting methods, and the analysis of these damping and cancellations effects
can be based on a general criterion, see for instance \cite[Chap.\,IV]{HuVe03}.
Here we will formulate such a criterion for multistep methods, which was
already used --\,in a slightly hidden form\,-- in \cite{Hu01} for a class of
adaptive implicit-explicit two-step methods.

We consider a stable error recursion in $\R^M$,
\begin{subequations}
\label{eq:ConvCrit}
\eq
\label{eq:ConvCrit1}
\eps_n = \sum_{i=1}^k R_i \, \eps_{n-i} \,+\, \delta_n \qquad
\mbox{for $n\ge k$}\,,
\eeq
with initial errors $\eps_0,\ldots,\eps_{k-1} = \Oh(\dt^k)$ and local
errors $\delta_n$ such that
\label{eq:ConvCrit2}
\eq
\begin{array}{c}
% \delta_n \,=\, \big( I - (R_1+\cdots+R_k) \big) \xi_n \,+\, \eta_n,
\displaystyle
\delta_n \,=\, \Big( I - \sum_{i=1}^k R_i\Big) \xi_n \,+\, \eta_n,
\\[6mm]
\xi_n = \Oh(\dt^k), \quad  \eta_n = \Oh(\dt^{k+1}),
\quad \xi_{n+1} - \xi_{n} = \Oh(\dt^{k+1}),
\end{array}
\eeq
\end{subequations}
uniformly for $n\ge k$.
Then $\eps_n = \Oh(\dt^k)$ uniformly for $t_n\in[0,T]$, that is, the
method is convergent of order $k$ on the interval $[0,T]$.

The proof of this statement is easy: defining $\eps_n^* = \eps_n - \xi_n$
for all $n$, where we set $\xi_n=\xi_k$ if $n<k$, it follows that
$$
\eps_n^* = \sum_{i=1}^k R_i \eps_{n-i}^* + \delta_n^*\,, \qquad
\delta_n^* = \sum_{i=1}^k R_i(\xi_{n-i}-\xi_n) + \eta_n \,,
$$
and for these transformed local errors we have $\delta_n^* = \Oh(\dt^{k+1})$.

To apply the convergence criterion (\ref{eq:ConvCrit2}), note that
the extrapolation coefficients in (\ref{eq:Lagr1}), (\ref{eq:Lagr2}) are
such that $\sum_{i=1}^k \widehat{c}_i = \sum_{i=1}^k \widecheck{c}_i = 1$. 
Consequently, if $\beta = \sum_{i=0}^k b_i$, then also
$\sum_{i=1}^k \widehat{b}_i = \sum_{i=1}^k \widecheck{b}_i = \beta$.
It will be tacitly assumed that the implicit method (\ref{eq:LMimpl}) is 
zero-stable and consistent, and then we will have $\sum_{i=1}^k a_i = 1$ 
and $\beta\neq0$.  Using (\ref{eq:R_i}), we therefore obtain
\eq
\label{eq:SumR_i}
\sum_{i=1}^k R_i \,=\, I + \beta P^{-1} Z
\eeq
with $Z = Z_0+Z_1+\cdots+Z_s$. Setting
\eq
\label{eq:zeta_n}
S_j \,=\, I - Q_1\cdots Q_{j-1} \,,
\qquad
\zeta_n = -\dt^k \mfrac{\beta_0}{\beta} \sum_{j=1}^s Z^{-1} S_j
\fhi_j^{(k-1)}(t_n) \,,
\eeq
it follows that the convergence criterion (\ref{eq:ConvCrit2}) can be
applied with $\xi_n = \zeta_n$ provided the terms $Z^{-1}S_j \upsilon(t)$
with $\upsilon(t) = \fhi_j^{(k-1)}(t), \fhi_j^{(k)}(t)$  are all bounded
uniformly in the stiffness. This is similar to the formulas obtained in
\cite{AHHP17} for a modified Douglas method, so we can repeat the main
arguments here.

Note that $S_1 = 0$, $S_2 = b_0 Z_1$, and for $2\le j\le s$ the following
expression is obtained:
\eq
\label{eq:S_j}
S_j = \sum_{m=1}^s \Big( (-1)^{m-1} \!\!\!
\sum_{1\le i_1<\cdots< i_m< j}
b_0^m Z_{i_1}Z_{i_2}\cdots Z_{i_m} \Big) \,.
\eeq
Consequently we will have $Z^{-1} S_j\upsilon(t) = \Oh(1)$ if all products
$Z^{-1}Z_{i_1}Z_{i_2}\cdots Z_{i_m}\upsilon(t) = \Oh(1)$ for
$1\le i_1<i_2<\cdots< i_m< j$. The essential condition for
a global error bound of order $k$ will therefore be:
\eq
\label{eq:ConvCond}
\begin{array}{r}
\dt^{m-1}  A^{-1} A_{i_1} A_{i_2} \cdots A_{i_m} \upsilon(t) = \Oh(1)
\qquad \mbox{for} \;\; \upsilon = \fhi^{(k-1)}_j, \fhi^{(k)}_j \,,\;
t \in [0,T] \;
\\
\mbox{and} \;
1 \le i_1 < \cdots < i_m < j \le s \,,
\end{array}
\eeq
with $A = A_0 + A_1 + \cdots + A_s$.
In summary, we have obtained the following convergence result.

\begin{Thm} \label{Thm:Conv}
Consider linear problems (\ref{eq:ODE}), (\ref{eq:ODElin}) on a finite
interval $[0,T]$ with $\fhi_j^{(k-1)}(t), \fhi_j^{(k)}(t) = \Oh(1)$
for $t\in[0,T]$, $j=2,3,\ldots,s$.
Assume the stability conditions (\ref{eq:BoundQinv}) and (\ref{eq:stability0})
are satisfied, and condition (\ref{eq:ConvCond}) holds.  Then method
(\ref{eq:SCLMmod}) will be convergent of order $k$ on the interval $[0,T]$.
\end{Thm}

If $s=2$, this result shows convergence with order $k$ under the condition
$A^{-1} A_1 \upsilon(t) = \Oh(1)$ for
$\upsilon = \fhi^{(k-1)}_2\!, \fhi^{(k)}_2$ with $t \in [0,T]$.
For $s=3$ we get the additional conditions $A^{-1} A_1 \upsilon(t)= \Oh(1)$,
$A^{-1} A_2 \upsilon(t) = \Oh(1)$ and
$\dt \, A^{-1} A_1 A_2 \upsilon(t) = \Oh(1)$
for $\upsilon = \fhi^{(k-1)}_3\!, \fhi^{(k)}_3$ with $t \in [0,T]$.
For a more detailed  discussion of these convergence conditions
for a 3D heat equation we refer to \cite{AHHP17,Hu02}, and in these
references it is also noted that non-singularity of $A$ is not essential.

\section{Stability for parabolic problems with mixed derivatives}
\label{Sect:Stab}

Parabolic equations with mixed derivatives arise, for example, in 
financial applications. 
As model problems to analyze stability for such applications we consider 
first in this section the 2D pure diffusion equation with mixed derivative 
and next the more general 2D advection-diffusion equation with mixed 
derivative.
These problems have previously been considered in the stability analysis 
of one-step splitting methods in for example 
\cite{CrSn88,HoMi13,HoWe07,HoWe09,McMi70}.

The 2D model pure diffusion equation is given by
\eq
\label{eq:finmod1}
u_t = d_{11} u_{x_1x_1} + d_{22} u_{x_2x_2} + (d_{12}+d_{21})\, u_{x_1x_2} 
\eeq
on the unit square $\Omega = [0,1]^2$ with periodic boundary conditions, 
and with coefficients $d_{ij}$ such that
\begin{subequations}
\label{eq:finmod2}
\eq
\label{eq:finmod2A}
D = (d_{ij}) \in \R^{2\times 2} \quad \mbox{is positive semi-definite}\,,
\eeq
\eq
\label{eq:finmod2B}
|d_{12}+d_{21}| \,\le\, 2\gamma\, \sqrt{d_{11} d_{22}} \,.
\eeq
\end{subequations}
The value of $\gamma \in [0,1]$ is a measure for the size of the
correlation factor of the two underlying stochastic processes in 
the financial model. 

For the spatial discretization standard central second-order finite 
differences are applied on uniform Cartesian grids with mesh-width 
$h>0$ in both directions. 
Considering an explicit treatment of the mixed derivative part 
followed by two implicit unidirectional corrections, we have a 
splitted linear ODE  system (\ref{eq:ODE}), (\ref{eq:Deco}) with 
$s=2$ and normal commuting matrices and real, scaled eigenvalues
\eq
\label{eq:finmod3}
z_j = -2 r\, d_{jj}\, (1-\cos\phi_j) \quad (j=1,2)\,, \qquad
z_0 = - r (d_{12}+d_{21}) \, \sin\phi_1 \, \sin\phi_2 \,,
\eeq
where $r = \dt/h^{2}$ and $\phi_j\in[0,2\pi]$; see e.g.\ \cite{HoWe09}. 
Using (\ref{eq:finmod3}), stability in the discrete $L_2$-norm follows 
for the ODE system obtained from (\ref{eq:finmod1}).
We are interested in unconditional stability of a given modified stabilizing 
correction splitting method, that is stability for  all $r>0$, and 
for a given method a sufficient condition on the parameter $\theta >0$ 
in function of $\gamma$ will be determined such that unconditional 
stability holds.
%To derive sufficient conditions, one can use $z_j \le 0$, $j=1,2,\ldots,s$,
%together with the following inequalities obtained in \cite{HoWe09,HoMi13}: 
%\eq
%\label{eq:finmod4}
%z \,\le\, 0 \,, \qquad
%|z_0| \,\le\, \gamma \sum_{i\neq j} \sqrt{z_i z_j} \,.
%\eeq
To derive these conditions, we use the properties
\eq
\label{eq:finmod4}
z_j \,\le\, 0 \quad (j=1,2)\,, \qquad
z\le 0\,, \qquad
|z_0| \,\le\, 2\gamma \, \sqrt{z_1 z_2} \,,
\eeq
with $z = z_0+z_1+z_2$, where the latter two inequalities were 
proved in \cite{HoMi13,HoWe09}.

Stability is relatively easy to study when $k=s=2$. 
Consider the polynomial
$\pi(\zeta) = \zeta^2 - r_1 \zeta - r_2$ with real coefficients $r_1,r_2$.
The two roots of $\pi$ both have modulus less than or equal to one iff
\eq
\label{eq:Schur}
|r_2| \,\le\, 1 \,, \qquad
|r_1| \,\le\, 1 - r_2 \,,
\eeq 
as can be seen, for example, by using the Schur criterion.
For the root condition, multiple roots of modulus one are to be excluded,
which happens if $r_1 = \pm 2$, $r_2 = -1$. In the following the stability 
criterion will be applied to the classes of two-step splitting methods that 
were introduced in Section~\ref{Sect:Exas}.
Recall that $p \,=\, \prod_{j=1}^s (1 - \theta z_j)$ and $\theta = b_0$.

\begin{Lem}\label{Lem:stab}
Let $s=2$ and let $\alpha, \beta_0, \beta$ be real numbers with $\alpha\ge -1$.
Then 
\eq
p + \alpha + \beta_0 z_0 + \beta (z_1+z_2) \ge 0
\eeq
whenever (\ref{eq:finmod4}) holds and
\eq
\label{eq:boundtheta}
\theta \ge \max \left\{ \beta\,,\,\frac{|\beta_0|\gamma + \beta}{1+\sqrt{1+\alpha}} \right\}.
\eeq 
\end{Lem}
\begin{proof}
Write $y_j = \sqrt{-\theta z_j}$ for $j=1,2$.
Then $|z_0| \,\le\, 2\displaystyle\frac{\gamma}{\theta} y_1 y_2$ and
$$
p + \alpha + \beta_0 z_0 + \beta (z_1+z_2) \ge 
(1+y_1^2)(1+y_2^2)  + \alpha - 2\kappa_0 y_1y_2 -\kappa(y_1^2+y_2^2)
$$
with $\kappa_0 = \displaystyle\frac{|\beta_0|\gamma}{\theta}$ and $\kappa = \displaystyle\frac{\beta}{\theta}$.
There holds
\begin{eqnarray*}
&& (1+y_1^2)(1+y_2^2) + \alpha - 2\kappa_0 y_1y_2 -\kappa(y_1^2+y_2^2) = \\
&&  1 + \alpha + y_1^2+y_2^2 + y_1^2y_2^2 - 2\kappa_0 y_1y_2 -\kappa(y_1^2+y_2^2) = \\
&& 1 + \alpha + (1-\kappa)(y_1-y_2)^2 + 2(1-\kappa)y_1y_2 + y_1^2y_2^2 - 2\kappa_0 y_1y_2 = \\
&& 1 + \alpha + (1-\kappa)(y_1-y_2)^2 + y_1^2y_2^2 + 2(1-\kappa_0-\kappa)y_1y_2 = \\
&& 1 + \alpha + (1-\kappa)(y_1-y_2)^2 + (y_1y_2 + 1-\kappa_0-\kappa)^2 - (1-\kappa_0-\kappa)^2.
\end{eqnarray*}
Using the latter two expressions it follows that $p + \alpha + \beta_0 z_0 + \beta (z_1+z_2)\ge 0$ whenever
$\kappa \le 1$ and $\kappa_0+\kappa \le 1+\sqrt{1+ \alpha}$.
Inserting $\kappa_0$, $\kappa$ yields the result of the lemma.
\end{proof}

\begin{Thm}\label{Thm:stab}
Consider (\ref{eq:finmod1}), (\ref{eq:finmod2}) with $s=2$ and periodic boundary condition.
Let (\ref{eq:ODE}), (\ref{eq:Deco}) be obtained by central second-order finite difference 
discretization and splitting as described above.
Then the three modified stabilizing correction methods (\ref{eq:SCLMmod}) given by $k=2$ 
and (\ref{eq:CNLF}), (\ref{eq:BDF2}), (\ref{eq:Adams2}), respectively, are unconditionally 
stable for the following parameter values $\theta$:
\noindent
\begin{itemize}
\item method (\ref{eq:CNLF}):~~~~$\theta = 1$,
\item method (\ref{eq:BDF2}):~~~~$\displaystyle \theta \geq\max\left\{\frac{1}{2}\,,\, \frac{\gamma+1}{2+2/\sqrt{3}} \right\}$\,,
\item method (\ref{eq:Adams2}):~~~~$\displaystyle \theta \geq\max\left\{\frac{1}{2}\,,\, \frac{\gamma+1}{3} \right\}$\,.
\end{itemize}
\end{Thm}
\noindent
We remark that in numerical experiments a smaller value $\theta$ is often seen to yield
smaller error constants.

\begin{proof}
(i) The first condition from (\ref{eq:Schur}) is equivalent to
$$
|a_2 + \widehat{b}_2 z_0 + \widecheck{b}_2 (z_1+z_2) + \mfrac{1}{b_0}(b_2-\widecheck{b}_2) (1 - p)| \le |p|.
$$
Since $b_2=\widecheck{b}_2$ for all three methods under consideration and $p\ge 1>0$, this holds iff
\begin{subequations}
\label{eq:cond1}
\eq
\label{eq:cond1A}
p + a_2 + \widehat{b}_2 z_0 + \widecheck{b}_2 (z_1+z_2) \ge 0,
\eeq
\eq
\label{eq:cond1B}
p - a_2 - \widehat{b}_2 z_0 - \widecheck{b}_2 (z_1+z_2) \ge 0.
\eeq
\end{subequations}
By Lemma \ref{Lem:stab}, the conditions (\ref{eq:cond1A}), (\ref{eq:cond1B}) are fulfilled if $|a_2|\le 1$ and
\begin{subequations}
\label{eq:cond1t}
\eq
\label{eq:cond1At}
\theta \ge \max \left\{ \widecheck{b}_2\,,\,\frac{|\widehat{b}_2|\gamma + \widecheck{b}_2}{1+\sqrt{1+a_2}} \right\},
\eeq
\eq
\label{eq:cond1Bt}
~~\theta \ge \max \left\{ -\widecheck{b}_2\,,\,\frac{|\widehat{b}_2|\gamma - \widecheck{b}_2}{1+\sqrt{1- a_2}} \right\}.
\eeq
\end{subequations}
It is readily verified that $|a_2|\le 1$ and (\ref{eq:cond1At}) are always satisfied for all three methods.
Next, (\ref{eq:cond1Bt}) is always satisfied for method (\ref{eq:CNLF}).
For method (\ref{eq:BDF2}), condition (\ref{eq:cond1Bt}) reads
$$
\theta \ge \max \left\{ \frac{2}{3}-\theta\,,\,\frac{\sfrac{2}{3}\gamma + \sfrac{2}{3}-\theta}{1+2/\sqrt{3}} \right\},
$$
which is equivalent to
\eq
\label{eq:BDF2stab1}
\theta \ge \max \left\{ \frac{1}{3} \,,\, \frac{\gamma + 1}{3+\sqrt{3}} \right\}.
\eeq
Similarly, for method (\ref{eq:Adams2}) condition (\ref{eq:cond1Bt}) reads
$$
\theta \ge \max \left\{ \frac{1}{2}-\theta\,,\,\frac{\sfrac{1}{2}\gamma + \sfrac{1}{2}-\theta}{2} \right\},
$$
which is equivalent to
\eq
\label{eq:Adams2stab1}
\theta \ge \max \left\{ \frac{1}{4} \,,\, \frac{\gamma + 1}{6} \right\}.
\eeq
\noindent
(ii) The second condition from (\ref{eq:Schur}) holds iff
$$
|a_1 + \widehat{b}_1 z_0 + \widecheck{b}_1 (z_1+z_2) + p - 1|
+ a_2 + \widehat{b}_2 z_0 + \widecheck{b}_2 (z_1+z_2)  \le p,
$$
where it has been used that $\widecheck{b}_1 = b_1 + \theta$.
This inequality is equivalent to
\begin{subequations}
\label{eq:cond2}
\eq
\label{eq:cond2A}
\phantom{2p +  }~ a_1 + a_2 -1 + (\widehat{b}_1 + \widehat{b}_2) z_0 + (\widecheck{b}_1 + \widecheck{b}_2) (z_1+z_2) \le 0,
\eeq
\eq
\label{eq:cond2B}
2p + a_1 - a_2 -1 + (\widehat{b}_1 - \widehat{b}_2) z_0 + (\widecheck{b}_1 - \widecheck{b}_2) (z_1+z_2) \ge 0.
\eeq
\end{subequations}
From $a_1 + a_2 = 1$, $\widehat{b}_1 + \widehat{b}_2 = \widecheck{b}_1 + \widecheck{b}_2 >0$ and $z=z_0+z_1+z_2\le 0$
it follows that (\ref{eq:cond2A}) is always satisfied.
Next, by Lemma \ref{Lem:stab}, condition (\ref{eq:cond2B}) is satisfied if $\alpha\ge -1$ and one has the lower bound 
(\ref{eq:boundtheta}) where
$$
\alpha  = \sfrac{1}{2}(a_1 - a_2 -1) = - a_2,\qquad
\beta_0 = \sfrac{1}{2}(\widehat{b}_1 - \widehat{b}_2),\qquad
\beta   = \sfrac{1}{2}(\widecheck{b}_1 - \widecheck{b}_2).
$$
For method (\ref{eq:CNLF}) this is readily seen to be true.
For method (\ref{eq:BDF2}) this holds if
$$
\theta \ge \max \left\{ 1-\theta\,,\,\frac{\gamma + 1-\theta}{1+2/\sqrt{3}} \right\},
$$
which is equivalent to
\eq
\label{eq:BDF2stab2}
\theta \ge \max \left\{ \frac{1}{2} \,,\, \frac{\gamma + 1}{2+2/\sqrt{3}} \right\}.
\eeq
Finally, for method (\ref{eq:Adams2}) this holds if
$$
\theta \ge \max \left\{ 1-\theta\,,\,\frac{\gamma + 1-\theta}{2} \right\},
$$
which is equivalent to
\eq
\label{eq:Adams2stab2}
\theta \ge \max \left\{ \frac{1}{2} \,,\, \frac{\gamma + 1}{3} \right\}.
\eeq
\noindent
(iii) Combining the results of part (i) and (ii), it follows that (\ref{eq:Schur}) 
is always fulfilled for method (\ref{eq:CNLF}) and it is fulfilled for method
(\ref{eq:BDF2}), respectively (\ref{eq:Adams2}), if the lower bound 
(\ref{eq:BDF2stab2}), respectively (\ref{eq:Adams2stab2}), holds.
It remains to show, for the root condition, that there are no multiple roots of 
modulus one.
To this purpose, suppose $r_2 = -1$, that is
$$
a_2 + \widehat{b}_2 z_0 + \widecheck{b}_2 (z_1+z_2) + p = 0.
$$
This yields
$$
1 + a_2 + \widehat{b}_2 z_0 + (\widecheck{b}_2-\theta) (z_1+z_2) + \theta^2 z_1z_2 = 0,
$$
and since $\widecheck{b}_2-\theta = \widehat{b}_2$,
$$
1 + a_2 + \widehat{b}_2 z + \theta^2 z_1z_2 = 0.
$$
Using that $z_1, z_2, z \le 0$, this yields a contradiction for each of the three
methods.
Thus the root condition is satisfied, which completes the proof of the theorem.
\end{proof}

In most financial applications the pertinent PDEs are of the advection-diffusion
kind and the lower bounds on $\theta$ derived above for unconditional stability 
may be too optimistic if advection is dominating.
Therefore, and with a view to the particular financial application in the next 
section, we consider also the 2D model advection-diffusion problem
\eq
\label{eq:finmod5}
u_t = c_1 u_{x_1} + c_2 u_{x_2} + d_{11} u_{x_1x_1} + d_{22} u_{x_2x_2} + 
(d_{12}+d_{21})\, u_{x_1x_2} 
\eeq
on the unit square $\Omega = [0,1]^2$ with periodic boundary conditions 
and coefficients $c_i$, $d_{ij}$ such that (\ref{eq:finmod2}) holds.
For the spatial discretization of (\ref{eq:finmod5}), again standard central 
second-order finite differences are applied on uniform Cartesian grids with 
mesh-width $h>0$ in both directions. 
The obtained semi-discrete ODE system (\ref{eq:ODE}) is splitted according 
to (\ref{eq:Deco}) with $s=2$ where $F_0$ represents the mixed derivative
part and $F_j$ represents all spatial derivatives in the $x_j$-direction
for $j=1,2$.
This leads to scaled eigenvalues
\eq
\label{eq:finmod6}
z_j = -2 r\, d_{jj}\, (1-\cos\phi_j) + \imi q\, c_j \sin\phi_j\,, \quad
z_0 = - r (d_{12}+d_{21}) \, \sin\phi_1 \, \sin\phi_2 \,,
\eeq
where $\imi = \sqrt{-1}$ and $q = \dt/h$ and $r$, $\phi_j$ are as before
($j=1,2$).
Using (\ref{eq:finmod2}), they are seen to satisfy, compare \cite{HoWe07},
\eq
\label{eq:finmod7}
\Real (z_1)\le 0, \quad \Real (z_2) \le 0 \quad {\rm and}\quad 
|z_0|\le 2\gamma\, \sqrt{\Real (z_1)\, \Real (z_2)}\,.
\eeq
We examine for the two-step modified stabilizing correction 
methods (\ref{eq:SCLMmod}) under consideration for which values of $\theta$ 
the root condition is satisfied whenever (\ref{eq:finmod7}) holds.
For complex coefficients $r_1, r_2$ the root condition is equivalent to 
\eq
\label{eq:Schur_complex}
|r_2| \,\le\, 1 \,, \qquad
|r_1+\bar{r}_1 r_2| \,\le\, 1 - |r_2|^2 \,,
\eeq 
and there are no multiple roots of modulus one.

For method (\ref{eq:CNLF}) the result is unfavourable.
Consider $z_0=0$ and $z_1=z_2=\imi y$ with $y\in \R$. 
Then the requirement $|r_2|\le 1$ becomes $|1+ 2\imi y| \le |1-\imi y|^2$,
which is easily seen to be violated whenever $0<|y|< \sqrt{2}$.
Hence, for method (\ref{eq:CNLF}), it already does not hold that the root 
condition is always fulfilled under (\ref{eq:finmod7}) if $\gamma=0$ 
(corresponding to no mixed derivative term).

For methods (\ref{eq:BDF2}), (\ref{eq:Adams2}) the result appears to 
be positive.
An analytical study for these methods of the root condition under 
(\ref{eq:finmod7}) is expected to be quite technical.
Therefore, we have conducted a numerical experiment to gain insight 
into the possible outcome.

\begin{figure}
\begin{center}
\includegraphics[width=0.75\textwidth]{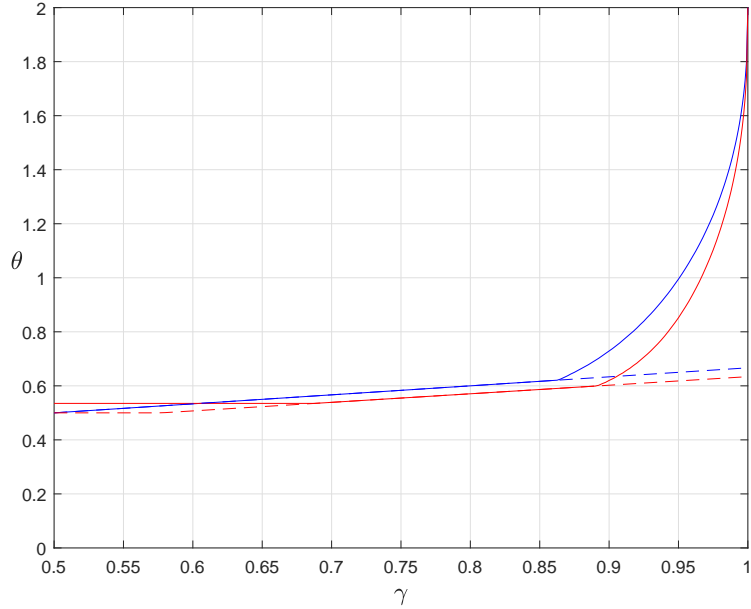}
\caption{\label{Fig:thetagamma} Solid curves: estimated lower 
stability bounds on $\theta$ versus $\gamma$ in the case of the 2D 
model advection-diffusion equation with mixed derivative. 
Dashed curves: analytical lower bounds given by Theorem~\ref{Thm:stab}
in the case of the 2D model pure diffusion equation with mixed derivative.
Red: method (\ref{eq:BDF2}). 
Blue: method (\ref{eq:Adams2}).}
\end{center}
\end{figure}

Let $w_{1,0}$ and $w_{i,j}$ for $i,j=1,2$ denote independent, uniformly 
distributed random numbers in $[0,1]$ and consider random triplets
$(z_0,z_1,z_2)$ and $(z_0,z_1,z_1)$ given by
\eq
z_j = -10^{1-5w_{1,j}}\pm \imi\,10^{1-5w_{2,j}}\,,
\quad
z_0 = (2w_{1,0}-1)\cdot2\gamma\,\sqrt{\Real (z_1)\, \Real (z_2)}
\eeq
for $j=1,2$.
Then (\ref{eq:finmod7}) holds and $z_0\in \R$.
Triplets with $z_1=z_2$ have been included as they are often found to 
yield the strongest requirement.
For each $\theta$ from a dense set of points in $[\sfrac{1}{2},2]$ we 
have estimated the maximal value $\gamma \in [0,1]$ (if any) such that 
the root condition is fulfilled whenever (\ref{eq:finmod7}) holds by 
testing the condition (\ref{eq:Schur_complex}) for two million random
triplets specified above.
The obtained numerical results for the two methods (\ref{eq:BDF2}), 
(\ref{eq:Adams2}) are shown in Figure \ref{Fig:thetagamma} as solid
red and blue curves, respectively, where $\theta$ has been displayed 
versus $\gamma \in [\sfrac{1}{2},1]$.
On the region $\gamma \in [0,\sfrac{1}{2}]$ the two curves are 
horizontal and this part is not shown.
The two curves represent estimated lower bounds on $\theta$, meaning
that unconditional stability - that is, without any restriction on $r$ 
or $q$ - is expected to hold in the application to (\ref{eq:finmod5}),  
(\ref{eq:finmod7}) whenever, for a given method and value $\gamma$, 
the value $\theta$ lies above the pertinent point on the curve.
%Based on Figure~\ref{Fig:thetagamma} we conjecture that for both
%methods (\ref{eq:BDF2}), (\ref{eq:Adams2}) there exists such a 
%lower bound on $\theta$ that is continuous and non-decreasing as
%a function of $\gamma$.
For comparison, the corresponding analytical lower bounds on $\theta$ from 
Theorem \ref{Thm:stab} for the 2D model pure diffusion problem have been 
included in the figure as dashed curves.
The estimated lower bounds in Figure \ref{Fig:thetagamma} for the 2D 
model advection-diffusion problem are close (or equal) to those up to 
$\gamma \approx 0.89$ and $\gamma \approx 0.86$ for methods (\ref{eq:BDF2}) 
and (\ref{eq:Adams2}), respectively.
Beyond this point, the lower bound for the 2D model advection-diffusion 
problem increases strongly for both methods, up to $\theta=2$ if
$\gamma = 1$.
In many financial applications, however, $\gamma$ is at most 0.9,
and this can be employed in the selection of a (smaller) value $\theta$.

For the stabilizing correction methods (\ref{eq:SCLM}), which use
only high-order extrapolation, a similar unconditional stability result 
in the case of the 2D model pure diffusion equation can be proved as 
Theorem \ref{Thm:stab} for the two-step modified stabilizing correction
methods (\ref{eq:SCLMmod}).
However, in the case of the 2D model advection-diffusion equation 
(\ref{eq:finmod5}) the stability results for the two-step methods 
(\ref{eq:SCLM}) appear to be much less favourable than those for 
the corresponding methods (\ref{eq:SCLMmod}).
Indeed, numerical experiments indicate that for the methods (\ref{eq:SCLM})
given by (\ref{eq:BDF2}) or (\ref{eq:Adams2}) (where $\widecheck{b}$ is now 
superfluous) unconditional stability is lacking already if $\gamma=0$.
In particular, for these methods, considering $z_0=0$ and $z_1=z_2=\imi y$ 
with $y\in \R$, numerical evidence suggests that for each $\theta\ge\frac{1}{2}$ 
there exists a value $y_0\in(0,2)$ such that there is instability whenever 
$|y|>y_0$.
A detailed stability study of the two-step methods (\ref{eq:SCLM}) for 
(\ref{eq:finmod5}) is beyond the scope of this paper, but it appears that 
at least an upper bound on $q = \dt/h$ is required to guarantee stability 
of these methods if advection is present.

We conclude with a brief discussion of the BDF3-type methods (\ref{eq:SCLM}) 
and (\ref{eq:SCLMmod}) given by (\ref{eq:BDF3}) when applied to the 2D model 
equation (\ref{eq:finmod5}).
The stability of these three-step methods has been investigated by a numerical 
study of the pertinent Schur criterion for complex polynomials of degree 3. 
For method (\ref{eq:SCLM}) given by (\ref{eq:BDF3}), we find that already
for the 2D model pure diffusion equation without mixed derivative ($\gamma=0$)
unconditional stability is lacking. 
More precisely, choosing $z_0=0$ and $z_1=z_2=y\in\R$ we find that there is 
instability whenever $y \lesssim -4.6$, which is clearly a negative result.
On the other hand, for method (\ref{eq:SCLMmod}) given by (\ref{eq:BDF3}), 
numerical evidence indicates that there is unconditional stability for the
2D model pure diffusion equation if $\gamma \lesssim 0.68$, which is a
favourable result.
For the general 2D model advection-diffusion equation, an upper bound on 
$q$ appears to be required however to guarantee stability. 
For example, taking $z_0=0$ and $z_1=z_2=\imi y$ with $y\in \R$, instability
is obtained whenever $|y|\gtrsim 2.7$.

\section{Test results for the Heston model}
\label{Sect:Heston}

To test the methods we consider a test set consisting of six parameter 
choices for the Heston model. These correspond to those in Haentjens \& 
in\,'t\,Hout \cite{HaHo12}.

\subsection{The Heston model}
The Heston model \cite{H93} for the fair values of European-style call 
options leads to a 2D time-dependent PDE of the form
\eq
\label{eq:Heston}
\frac{\partial u}{\partial t} = 
\sfrac{1}{2} s^2 v \frac{\partial^2 u}{\partial^2 s} +
\rho\sigma s v \frac{\partial^2 u}{\partial s \partial v} +
\sfrac{1}{2} \sigma^2 v \frac{\partial^2 u}{\partial^2 v} +
(r_{\! d} - r_{\! f}) s \frac{\partial u}{\partial s} +
\kappa (\eta - v) \frac{\partial u}{\partial v} -
r_{\!d} \, u 
\eeq
with independent variables $t\in (0,T]$ and $s,v > 0$, and initial 
condition 
\eq
\label{eq:HestonIC}
u(s, v, 0) \,=\, \max(0, s - K) \qquad (s\ge 0, v\ge 0).
\eeq
Here $T>0$ and $K>0$ are the given maturity date and strike price
of the option.
The parameter $\kappa>0$ is the mean-reversion rate, $\eta>0$ is 
the long-term mean, $\sigma>0$ is the volatility-of-variance, 
$\rho \in [-1,1]$ is the correlation between the two underlying
Brownian motions, and $r_d$, $r_f$ denote the domestic and foreign 
interest rates, respectively.
For feasibility of the numerical solution, the spatial domain is
truncated to a bounded set $[0,S_{\max}]\times [0,V_{\max}]$
with fixed values $S_{\max}$, $V_{\max}$ taken sufficiently large.
The following boundary conditions are imposed,
\noindent
\begin{subeqnarray}\label{eq:HestonBC}
\phantom{\frac{\partial}{\partial}}u(s,v,t)~~=&0
\quad &{\rm whenever}~~s=0\,,\\
\frac{\partial u}{\partial s}(s,v,t)~~=&e^{-r_{\! f} t}
\quad &{\rm whenever}~~s=S_{\max}\,,\\
\phantom{\frac{\partial u}{\partial s}}u(s,v,t)~~~=&se^{-r_{\! f} t}
\quad &{\rm whenever}~~v=V_{\max}\,.
\end{subeqnarray}
Further, at the $v=0$ boundary the PDE (\ref{eq:Heston}) is fulfilled, 
see \cite{ET11}.

The spatial discretization of the initial-boundary value problem for
(\ref{eq:Heston}) is performed using second-order finite differences 
on a smooth, nonuniform, Cartesian grid in the $(s,v)$-domain similar 
to that in \cite{HaHo12,HoFo10}.
The spatial grid has relatively many points in the neighbourhood of
the location $(s,v)=(K,0)$, which has been done both for financial
and numerical reasons.
We note that at the boundary $v=0$ the derivative $\partial u/\partial v$ 
is approximated using a second-order forward finite difference formula.
All other derivative terms in the $v$-direction vanish at $v=0$.
Cell averaging is applied to define the initial vector obtained from 
the nonsmooth initial condition (\ref{eq:HestonIC}) near the strike, 
see e.g.\ \cite{TR00}.
The resulting semi-discrete ODE system (\ref{eq:ODE}) is splitted according 
to (\ref{eq:Deco}) with $s=2$ where $F_0$ represents the mixed derivative
part and $F_1$, respectively $F_2$, represents all spatial derivatives in 
the $s$-direction, respectively $v$-direction, see e.g. \cite{HaHo12,HoFo10}.
For the subsequent numerical tests we choose the six cases of parameter 
sets for the Heston model listed in Table~\ref{Cases}. These correspond 
to those from \cite{HaHo12}.
 
\begin{table}[h]
\begin{center}
\begin{tabular}{|c|l|l|l|l|l|l|}
        \hline
        & Case A & Case B & Case C & Case D & Case E & Case F\\
        \hline \hline
        $\kappa$    & 3    & 0.6067  & 2.5    & 0.5   & 0.3   & 1\\
        $\eta$      & 0.12 & 0.0707  & 0.06   & 0.04  & 0.04  & 0.09\\
        $\sigma$    & 0.04 & 0.2928  & 0.5    & 1     & 0.9   & 1\\
        $\rho$      & 0.6  & -0.7571 & -0.1   & -0.9  & -0.5  & -0.3\\
        $r_d$       & 0.01 & 0.03    & 0.0507 & 0     & 0     & 0\\
        $r_f$       & 0.04 & 0       & 0.0469 & 0     & 0     & 0\\
        $T$         & 1    & 3       & 0.25   & 10    & 15    & 5\\
        $K$         & 100  & 100     & 100    & 100   & 100   & 100\\
        \hline
\end{tabular}
\end{center}
\caption{Parameter sets for the Heston model.}
\label{Cases}
\end{table}
\noindent
Cases A, B, C have previously been considered in \cite{HoFo10} and stem 
from \cite{B05,SST04,WAW02}, respectively.
Here the so-called Feller condition $2\kappa\eta > \sigma^2$ always holds.
A special feature of Case A is that $\sigma$ is close to zero, which implies 
that the PDE (\ref{eq:Heston}) is advection dominated in the $v$-direction.
Cases D, E, F were proposed in \cite{A08} as challenging test sets for 
practical applications.
In these three cases the maturity times are large and the Feller condition
is violated.

\subsection{One-step stabilizing correction methods}

The multistep methods will be compared with several well-known one-step
methods for problems (\ref{eq:ODE}), (\ref{eq:Deco}). The Douglas method
was already briefly introduced in Section~\ref{Sect:Exas}. Written out 
in full, the method reads
\eq
\label{eq:Do}
\setlength{\arraycolsep}{1mm}
\left\{
\begin{array}{ccl}
v_{0} &=& u_{n-1} + \dt\,F(t_{n-1}, u_{n-1}) \,,
\\[2mm]
v_{j} &=& v_{j-1}
+ \theta\dt \big(F_j(t_{n}, v_{j}) - F_j(t_{n-1}, u_{n-1})\big)
\qquad (j = 1,2,\ldots,s) \,,
\\[2mm]
u_{n} &=& v_{s} \,,
\end{array}
\right.
\eeq
with parameter $\theta \ge {1\over2}$.  Even if $\theta = {1\over2}$, 
the order is only one, due to the treatment of the explicit term $F_0$ 
in an Euler fashion. We note that the modification of this method that
was recently presented in \cite{AHHP17} is not sufficiently stable for 
the Heston model with explicit treatment of the cross-derivatives.

An extension of the Douglas method, due to in\,'t\,Hout \& Welfert 
\cite{HoWe09}, is given by
\eq
\label{eq:MCS}
\setlength{\arraycolsep}{1mm}
\left\{
\begin{array}{ccl}
v^*_{0} &=& u_{n-1} + \dt\,F(t_{n-1}, u_{n-1}) \,,
\\[2mm] 
v^*_{j} &=& v^*_{j-1}
+ \theta\dt \big(F_j(t_{n}, v^*_{j}) - F_j(t_{n-1}, u_{n-1})\big)
\qquad (j = 1,2,\ldots,s) \,,
\\[2mm]
v_{0} &=& v^*_{0} +
\sfrac{1}{2}\dt \big(F_0(t_{n}, v^*_{s}) - F_0(t_{n-1}, u_{n-1})\big)
\\[1mm]
& & \; + \, \big(\sfrac{1}{2}-\theta\big) \dt
\displaystyle
\sum_{j=1}^s \big(F_j(t_n, v^*_s) - F_j(t_{n-1}, u_{n-1}) \big) \,,
\\[2mm]
v_{j} &=& v_{j-1}
+ \theta\dt \big(F_j(t_{n}, v_{j}) - F_j(t_{n-1}, u_{n-1})\big)
\qquad (j = 1,2,\ldots,s) \,,
\\[2mm]
u_{n} &=& v_{s} \,.
\end{array}
\right.
\eeq
If $\theta = \frac{1}{2}$ this is the method of Craig \& Sneyd
\cite{CrSn88}.  Taking $\theta \in [\frac{1}{4},\frac{1}{2})$ often gives
better accuracy; in the numerical tests we will consider $\theta=\frac{1}{3}$.
For any choice of~$\theta$, method (\ref{eq:MCS}) is of order two in the 
ODE sense. Convergence results for PDEs have been derived in \cite{HoWy14}.
Stability results for (\ref{eq:MCS}) applied to the model problems 
(\ref{eq:finmod1}), (\ref{eq:finmod5}) can be found in 
\cite{HoMi11,HoMi13,HoWe09}, for example.

\subsection{Results for the Heston model} 

In the numerical tests we will compare the multistep methods with the 
following one-step methods:
\begin{center}
\begin{tabular}{l}
{Do} : the Douglas method (\ref{eq:Do}) with $\theta = \frac{1}{2}$,
\\[1mm]
{CS} : the Craig-Sneyd method, given by (\ref{eq:MCS}) with 
$\theta = \frac{1}{2}$,
\\[1mm]
{MCS} : the modified Craig-Sneyd method (\ref{eq:MCS}) with 
$\theta = \frac{1}{3}$.
\end{tabular}
\end{center}
Together with these well-known one-step methods we consider the following
two-step methods (\ref{eq:SCLMmod}) with modified stabilizing corrections:
\begin{center}
\begin{tabular}{l}
{SC2A} : the Adams2-type method (\ref{eq:Adams2}) with
$\theta = \frac{3}{4}$,
\rule{11mm}{0mm}
\\[1mm]
{SC2B} : the BDF2-type method (\ref{eq:BDF2}) with
$\theta = \frac{2}{3}$,
\\[1mm]
{SC2C} : the CNLF-type method (\ref{eq:CNLF}) with
$\theta = 1$.
\end{tabular}
\end{center}
In these Adams2- and BDF2-type methods the parameter value $\theta$ was chosen
so as to give a reasonable compromise between stability properties and error
constants.  The same holds for the modified Craig-Sneyd method.

\begin{figure}[ht!]
\setlength{\unitlength}{1cm}
\begin{center}
\begin{picture}(6,4.9)
\includegraphics[width=6cm]{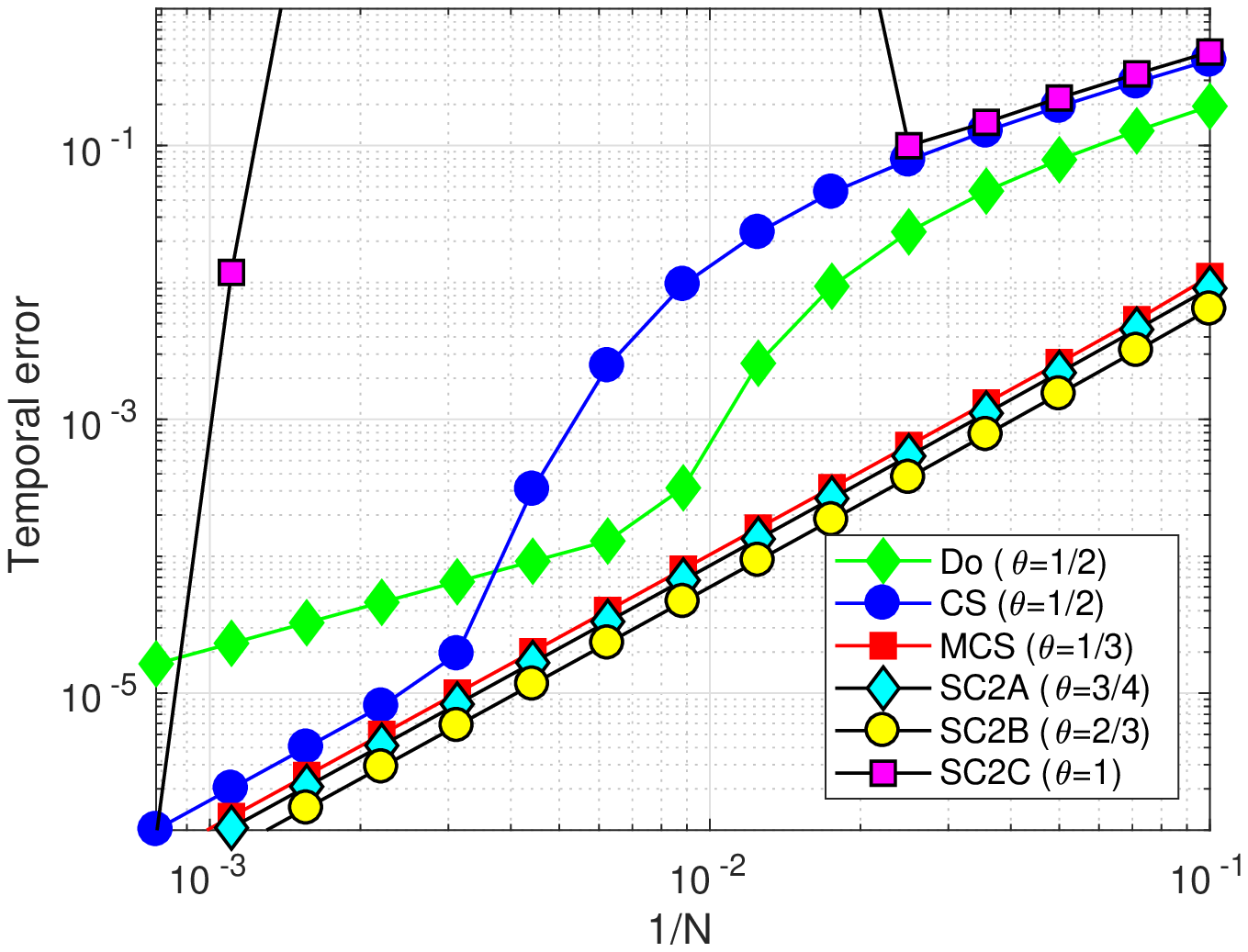}
\put(-5,4.2){\footnotesize Case A}
\end{picture}
\hspace{0.5cm}
\begin{picture}(6,4.9)
\includegraphics[width=6cm]{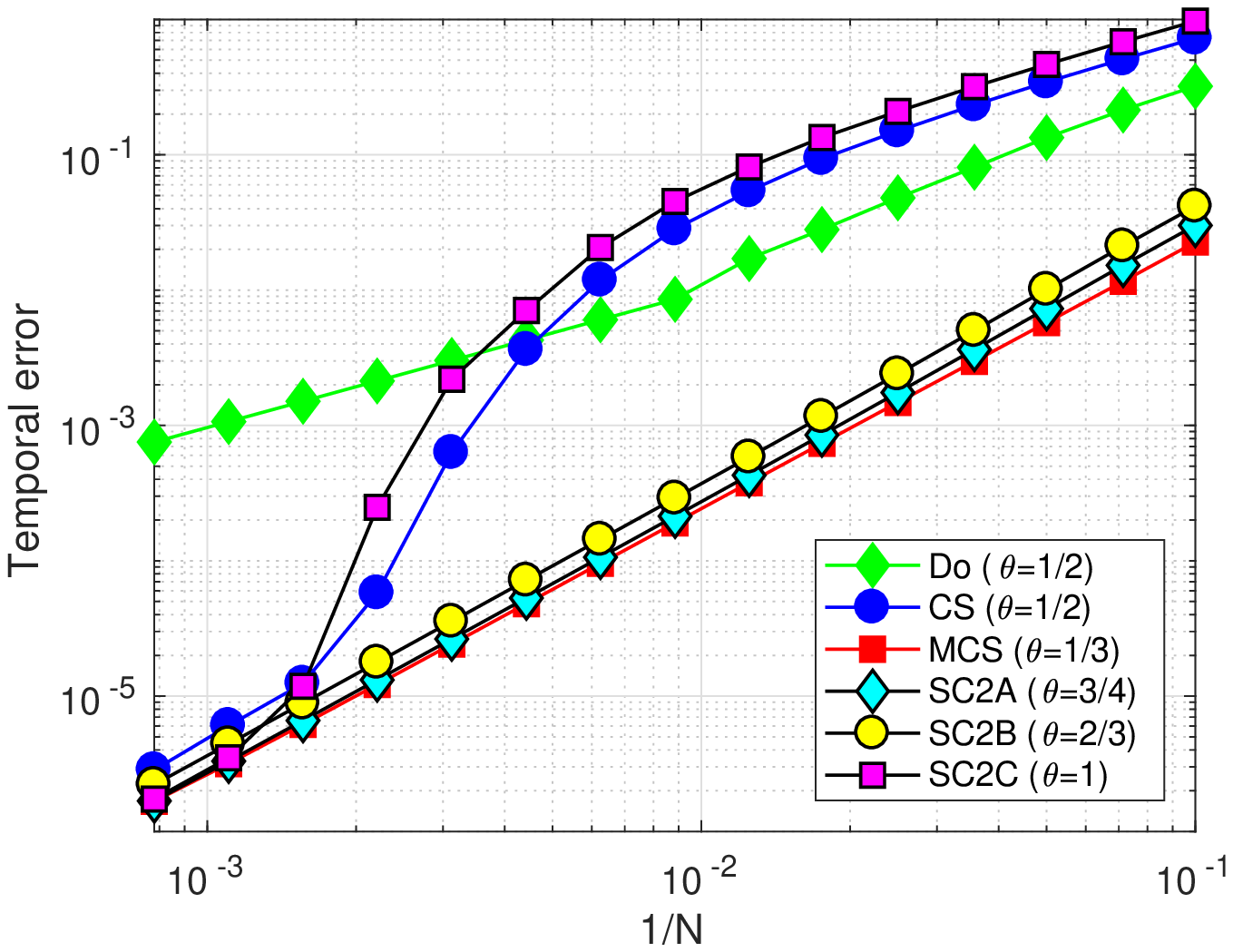}
\put(-5,4.2){\footnotesize Case B}
\end{picture}
\\[2mm]
\begin{picture}(6,4.9)
\includegraphics[width=6cm]{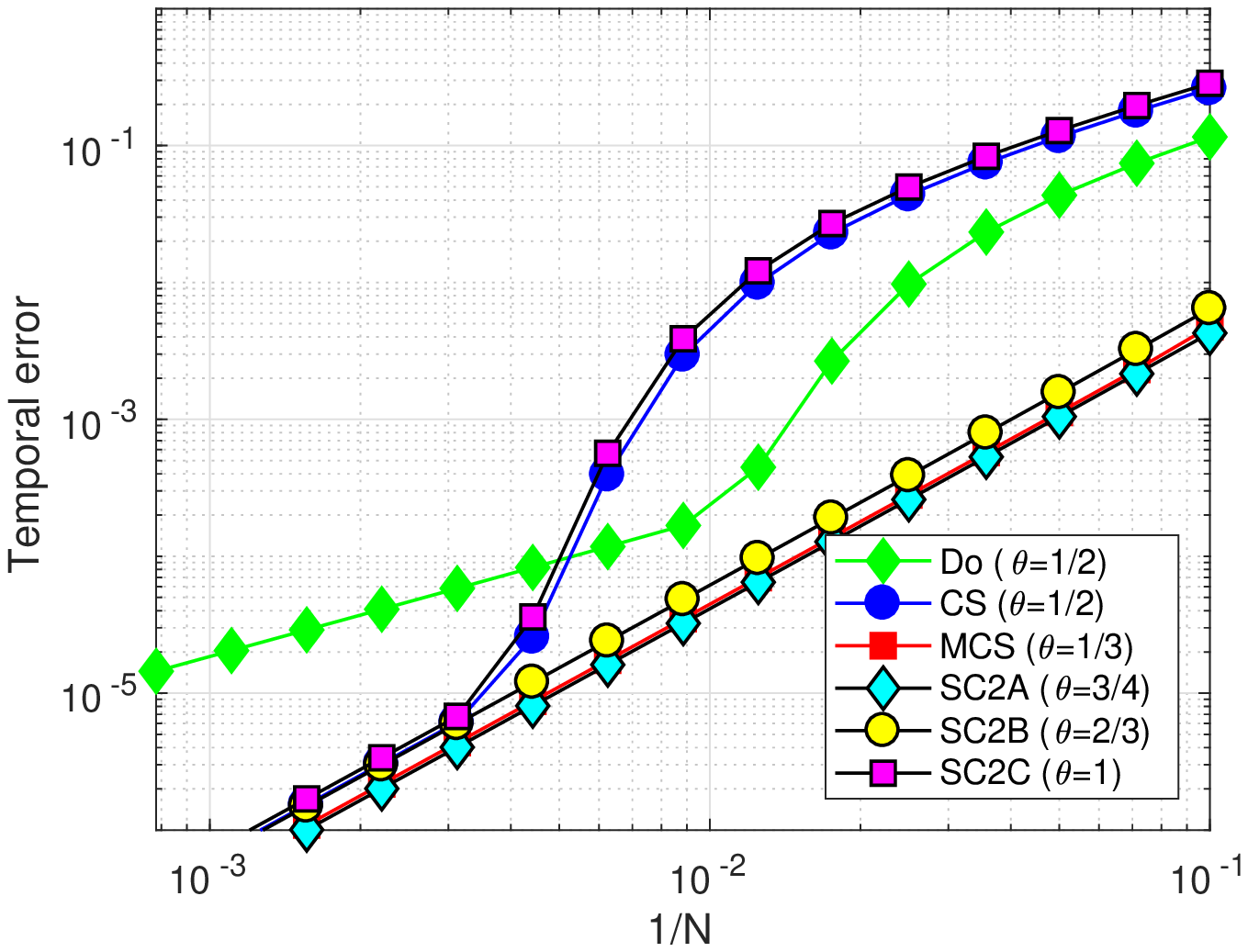}
\put(-5,4.2){\footnotesize Case C}
\end{picture}
\hspace{0.5cm}
\begin{picture}(6,4.9)
\includegraphics[width=6cm]{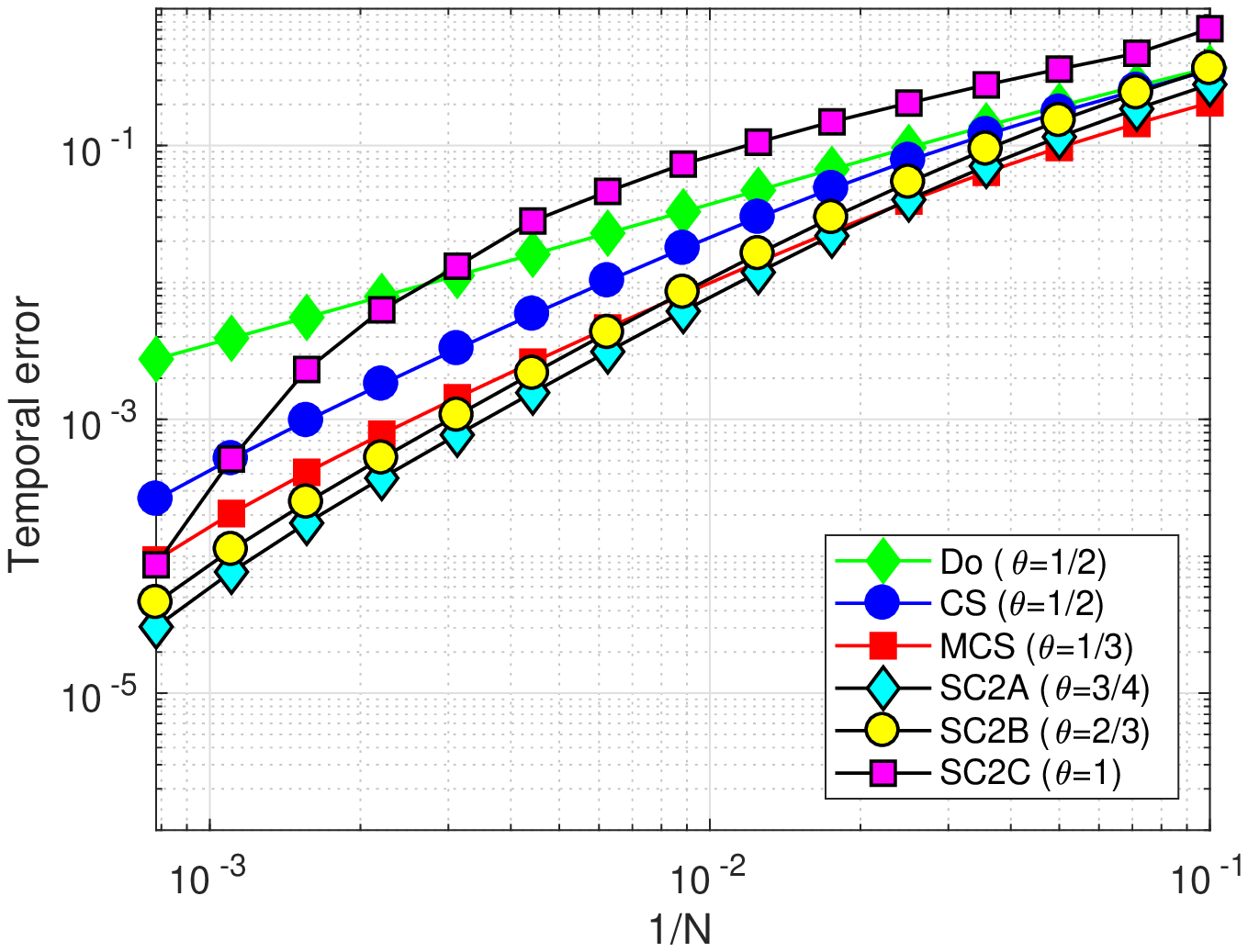}
\put(-5,4.2){\footnotesize Case D}
\end{picture}
\\[2mm]
\begin{picture}(6,4.9)
\includegraphics[width=6cm]{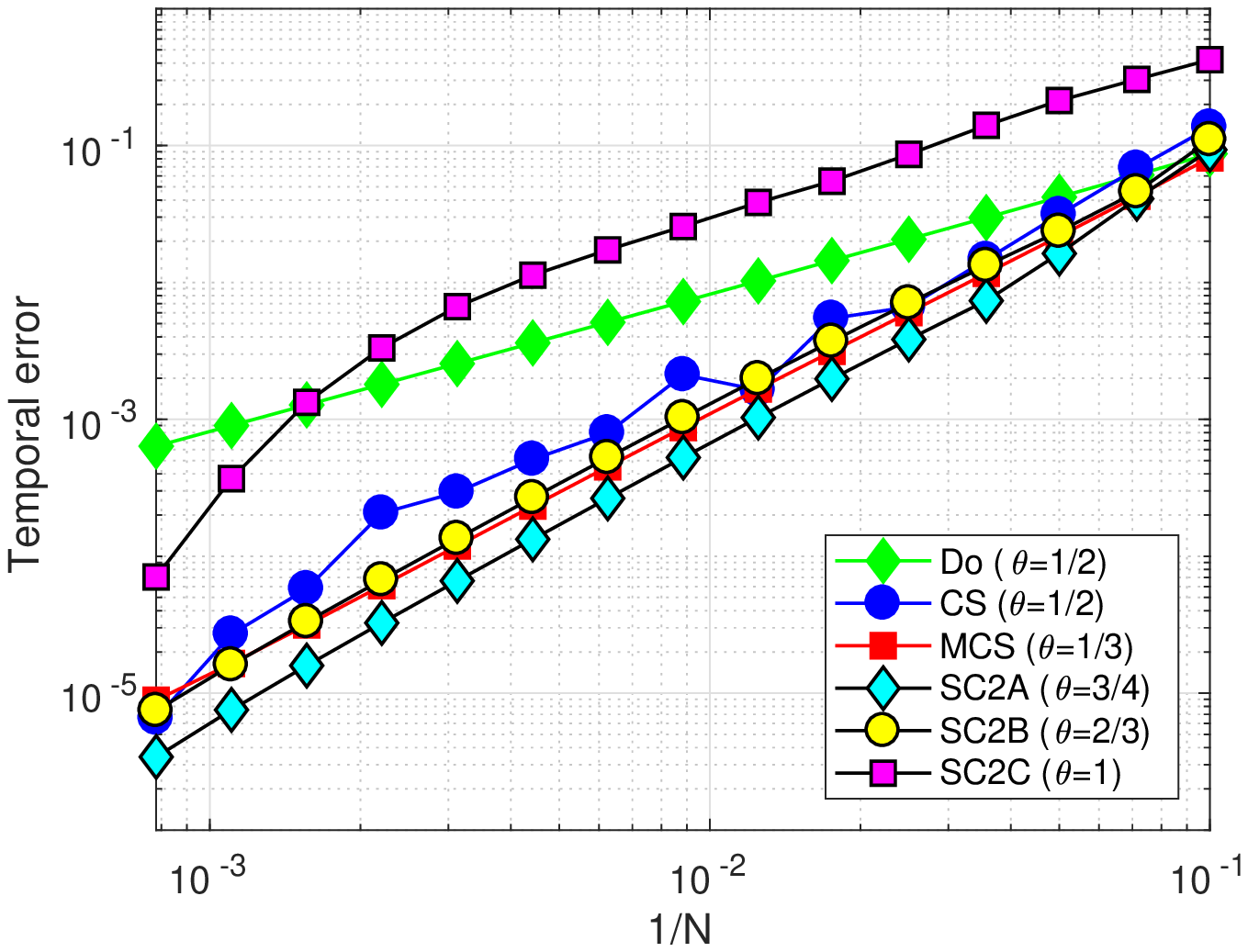}
\put(-5,4.2){\footnotesize Case E}
\end{picture}
\hspace{0.5cm}
\begin{picture}(6,4.9)
\includegraphics[width=6cm]{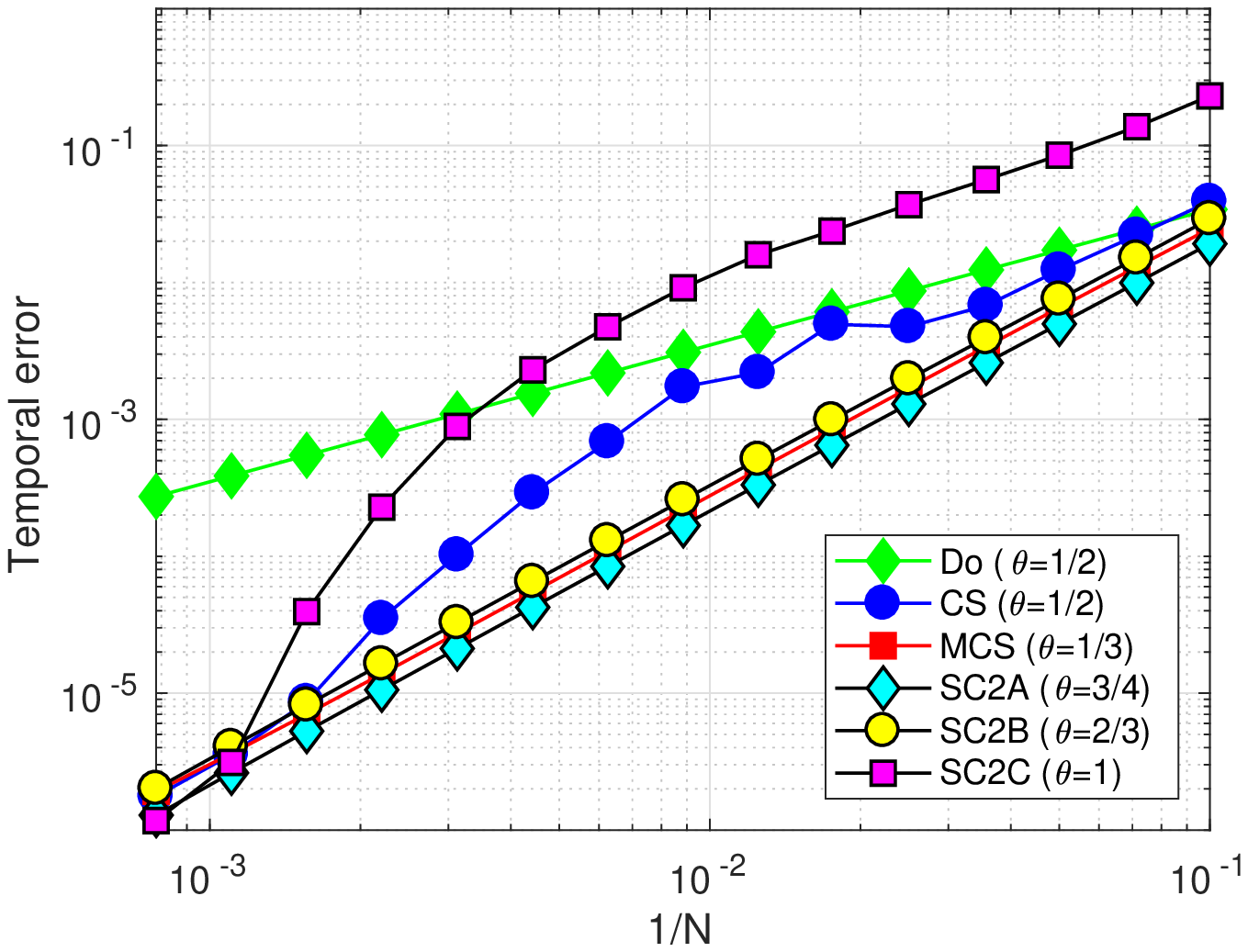}
\put(-5,4.2){\footnotesize Case F}
\end{picture}
\caption{ \small  \label{Fig:TestHeston}
Global errors at $t=T$ versus $1/N$ for the test cases A--F.
Methods CS and MCS with step-size $\dt = T/N$,
Douglas method and the multistep methods with step-size $\dt = T/[2 N]$.
}
\end{center}
\end{figure}

Per step, the methods CS and MCS are twice as expensive as the others, and
therefore we will use these methods with a step-size two times larger than
for the Douglas method and the multistep methods. 
For the two-step methods, the first approximation $u_1$ was computed
with the Douglas method with $\theta=1$, which seems a natural starting
method for the stabilizing correction two-step methods.

In Figure~\ref{Fig:TestHeston} the results are found for the six cases
of Heston parameter sets given by Table~\ref{Cases}.
In these plots, the global errors in the maximum norm are plotted as a
function of $1/N$, where $T/N$ is the step-size used for the CS and MCS 
methods.
Here a $400\times 200$ spatial grid has been taken and the global error 
is considered for $t=T$ on a region of financial interest given by 
$\frac{1}{2}K < s < \frac{3}{2}K$ and $0<v<1$.
Note that the global error does not contain the error due to spatial 
discretization.

From Figure~\ref{Fig:TestHeston}
it is seen that the two-step methods with stabilizing corrections based 
on the implicit Adams method (SC2A) and implicit BDF2 method (SC2B) are 
competitive with the modified Craig-Sneyd method (MCS).

The CNLF method (SC2C) behaves quite poorly in all six cases, with
large errors for moderate step-sizes and an irregular error behaviour,
probably due to instability.
Indeed, in Section \ref{Sect:Stab} it was noticed that for the 2D model 
advection-diffusion problem already a small amount of advection can 
render this method unstable.

\begin{Rem}[Smoothing steps]  \rm
The accuracy of the one-step methods in these tests can be somewhat
improved by performing two (non-splitted) backward Euler sub-steps
with step-size $\frac{1}{2}\dt$ to compute the first approximation
$u_1 \approx u(\dt)$. These are the so-called Rannacher smoothing steps
\cite{R84}.
The positive effect was most pronounced for the methods Do and CS, 
but even with such smoothing steps these two methods are not competitive 
with the best methods in these tests (MCS, SC2A and SC2B).
For a simple comparison of the methods, with comparable work for all
methods, such (non-splitted) backward Euler smoothing steps were not
used in the tests. Moreover, at the moment, it is not clear how a 
proper smoothing procedure should be constructed for the multistep methods.
\hfill $\Diamond$
\end{Rem}

\begin{Rem}[Two-step methods (\ref{eq:SCLM})]  \rm
The lower order extrapolation in (\ref{eq:SCLMmod}) was introduced to
improve stability of the schemes. In the above tests this enhanced
stability was found to be necessary, in particular for the Adams2-type 
scheme. The BDF2-type method showed a more stable behaviour but also
that method failed for the advection dominated PDE given by Case~A. 

For step-sizes that are sufficiently small for having stability, 
say $\dt\le\tau_{\rm stab}$,
the accuracy of these two-step methods (\ref{eq:SCLM}) was in general 
slightly better than for the corresponding methods (\ref{eq:SCLMmod}).
However, the stability threshold $\tau_{\rm stab}$ depends on the problem 
and on the spatial mesh-width.

Since robustness is an important quality in option valuation applications,
the modified methods (\ref{eq:SCLMmod}) are, in our opinion, preferable
over the schemes (\ref{eq:SCLM}).
\hfill $\Diamond$
\end{Rem}

\begin{Rem}[Three-step methods (\ref{eq:SCLMmod})]  \rm
Tests were also performed with the stabilizing corrections BDF3-type 
scheme (\ref{eq:SCLMmod}) with coefficients (\ref{eq:BDF3}). Starting 
values $u_1$ and $u_2$ were computed with the MCS method. For smaller 
step-sizes this SC3B method was seen to give higher accuracy than the 
two-step methods SC2A and SC2B, but instabilities were again observed 
for larger step-sizes, making this method not suitable for the Heston 
problem.  (Needless to say, the BDF3-type scheme (\ref{eq:SCLM}) with 
high-order extrapolation turned out to be very unstable.)
\hfill $\Diamond$
\end{Rem}

\section{Concluding remarks}
\label{Sect:Concl}

Among the two-step methods with stabilizing corrections (\ref{eq:SCLMmod}) 
considered in this paper the behaviour of the methods based on BDF2 and 
Adams2 was satisfactory, and these methods appear to be competitive with 
the well-established modified Craig-Sneyd method for the Heston problem.
The BDF3-type scheme (\ref{eq:SCLMmod}) may be suited for other applications, 
such as reaction-diffusion problems, which is left for future research.

The stabilizing correction methods studied in this paper can be
viewed as generalizations of IMEX linear multistep methods. Such
IMEX multistep methods have been examined for 1D option valuation
models with jumps in \cite{SaTo14}. Based on accuracy, the authors
had a slight preference for the CNAB method, i.e.\ (\ref{eq:Adams2})
with $\theta=\frac{1}{2}$, over the BDF2 method. Subsequently, this
CNAB method was applied to 2D models in \cite{SaToLy14}, but without
dimension splitting. It is part of our research plans to examine
the behaviour of the stabilizing correction multistep schemes to
models with jumps together with dimension splitting.

\end{document}